\documentclass[11pt]{amsart}
\usepackage{amsxtra}
\usepackage{amssymb}
\theoremstyle{plain}

\newcommand{\cleqn}{\setcounter{equation}{0}}
\newcommand{\clth}{\setcounter{theorem}{0}}
\newcommand {\sectionnew}[1]{\section{#1}\cleqn\clth}
\newcommand{\nn}{\hfill\nonumber}
\newtheorem{theorem}{Theorem}[section]
\newtheorem{lemma}[theorem]{Lemma}
\newtheorem{definition-theorem}[theorem]{Definition-Theorem}
\newtheorem{proposition}[theorem]{Proposition}
\newtheorem{corollary}[theorem]{Corollary}
\newtheorem{definition}[theorem]{Definition}
\newtheorem{example}[theorem]{Example}
\newtheorem{remark}[theorem]{Remark}
\newtheorem{conjecture}[theorem]{Conjecture}
\newtheorem{notation}[theorem]{Notation}
\newcommand \bth[1] { \begin{theorem}\label{t#1} }
\newcommand \ble[1] { \begin{lemma}\label{l#1} }

\newcommand \bpr[1] { \begin{proposition}\label{p#1} }
\newcommand \bco[1] { \begin{corollary}\label{c#1} }
\newcommand \bde[1] { \begin{definition}\label{d#1}\rm }
\newcommand \bex[1] { \begin{example}\label{e#1}\rm }
\newcommand \bre[1] { \begin{remark}\label{r#1}\rm }
\newcommand \bcj[1] { \begin{conjecture}\label{j#1}\rm }

\newcommand \bnota[1] { \begin{notation}\label{n#1}\rm }
\renewcommand {\eth} { \end{theorem} }
\newcommand {\ele} { \end{lemma} }

\newcommand {\epr} { \end{proposition} }
\newcommand {\eco} { \end{corollary} }
\newcommand {\ede} { \end{definition} }
\newcommand {\eex} { \end{example} }
\newcommand {\ere} { \end{remark} }
\newcommand {\ecj} { \end{conjecture} }

\newcommand {\enota} { \end{notation} }
\newcommand \thref[1]{Theorem \ref{t#1}}
\newcommand \leref[1]{Lemma \ref{l#1}}
\newcommand \prref[1]{Proposition \ref{p#1}}
\newcommand \coref[1]{Corollary \ref{c#1}}

\newcommand \deref[1]{Definition \ref{d#1}}
\newcommand \exref[1]{Example \ref{e#1}}

\newcommand \Cset {{\mathbb C}}        
\newcommand \KK {{\mathbb K}}
\newcommand \Zset {{\mathbb Z}}
\newcommand \calA {{\mathcal{A}}}
\newcommand \NN {{\mathcal{N}}}
\newcommand \NGZ {{\mathcal{NGT}}}
\newcommand \calT {{\mathcal{T}}}

\newcommand \OO {{\mathcal{O}}}

\newcommand \UU {{\mathcal{U}}}
\newcommand \ZZ {{\mathcal{Z}}}
\newcommand \qb {{\bf{q}}}
\newcommand \de {\delta}
\newcommand \al {\alpha}

\newcommand \la {\lambda}
\newcommand \La {\Lambda}

\newcommand \sig {\sigma}
\newcommand \ci  {\circ}
\newcommand \rcor {\rangle}
\newcommand \lcor {\langle}
\newcommand \ol {\overline}

\newcommand \id { {\mathrm{id}} }
\DeclareMathOperator \rank { {\mathrm{rank}} }
\newcommand \g  {\mathfrak{g}}   

\DeclareMathOperator \Aut { {\mathrm{Aut}} }

\DeclareMathOperator \lt  { {\mathrm{lt}} }
\DeclareMathOperator \lc  { {\mathrm{lc}} }

\DeclareMathOperator \Fract { {\mathrm{Fract}} }

\newcommand\kx{\KK^*}

\newcommand\HH{{\mathcal{H}}}
\newcommand\xh{X(\HH)}
\DeclareMathOperator \Spec {Spec}

\DeclareMathOperator \chr {char}
\newcommand \Znn {\Zset_{\ge 0}}
\DeclareMathOperator \ltq  { {\mathrm{lt}_{\qb}} }
\newcommand \xbar {\ol{x}}
\newcommand \Hmax {\HH_{\max}}

\begin{document}
\title[From quantum Ore extensions to quantum tori]
{From quantum Ore extensions to quantum tori via noncommutative UFD{\tiny{s}}}
\author[K. R. Goodearl]{K. R. Goodearl}
\address{
Department of Mathematics \\
University of California\\
Santa Barbara, CA 93106 \\
U.S.A.
}
\email{goodearl@math.ucsb.edu}
\author[M. T. Yakimov]{M. T. Yakimov}
\thanks{The research of K.R.G. was partially supported by NSF grant DMS-0800948, 
and that of M.T.Y. by NSF grant DMS-1303038 and Louisiana Board of Regents grant Pfund-403}
\address{
Department of Mathematics \\
Louisiana State University \\
Baton Rouge, LA 70803 \\
U.S.A.
}
\email{yakimov@math.lsu.edu}
\date{}
\begin{abstract} All iterated skew polynomial extensions 
arising from quantized universal enveloping algebras of Kac--Moody
algebras are special examples of a very large, axiomatically defined 
class of algebras, called CGL extensions. For the purposes of  
constructing initial clusters for quantum cluster algebra structures on an algebra $R$, 
and classification of the
automorphisms of $R$, one needs embeddings of $R$ into quantum tori $\calT$
which have the property that $R$ contains the corresponding quantum affine space 
algebra $\calA$. We explicitly construct such an embedding 
$\calA \subseteq R \subset \calT$ for each CGL extension $R$ using the methods of noncommutative 
noetherian unique factorization domains and running a Gelfand--Tsetlin
type procedure with normal, instead of central elements.
Along the way we classify the homogeneous prime elements of all CGL extensions and we prove that each CGL extension $R$ has an associated maximal torus which covers the automorphisms of $R$ corresponding to all normal elements.
For symmetric CGL extensions, we describe the relationship between our quantum 
affine space algebra $\calA$ and Cauchon's quantum affine space algebra generated by elements obtained via deleting derivations.    
\end{abstract}

\dedicatory{Dedicated to the memory of Andrei Zelevinsky}

\maketitle

\sectionnew{Introduction}
\label{intro}
Let $R$ be a noncommutative (right, say) Ore domain over a field $\KK$. 
The first step of the construction 
of a quantum cluster algebra structure on $R$ is the construction of 
an initial cluster. This amounts to the construction of a chain of embeddings
\begin{equation}
\label{eemb}
\calA \subseteq R \subset \calT \subset \Fract(R)
\end{equation}
where $\calA$ is a quantum affine space algebra, $\calT$ is the corresponding 
quantum torus (see Subsection \ref{3.2} for details), and 
$\Fract(R)$ is the Ore division ring of $R$. Embeddings of the form 
\eqref{eemb} also play an important role in classifying $\Aut(R)$ 
and proving rigidity results for $R$ in a general scheme recently 
developed by the second author \cite{Y2,Y3}.
Assume that $R$ is a $\Znn$-graded algebra and $\calA$ and $\calT$ 
are equipped with $\Zset$-gradings in which their generators have positive degrees 
and such that the first two inclusions in \eqref{eemb} are graded. 
We call an automorphism $\varphi$ of $R$ \emph{unipotent} if 
$\varphi(a) - a \in R^{m+1}+ R^{m+2} + \cdots$ for all $a \in R^m$, 
$m \in \Znn$. By \cite[Proposition 3.3]{Y3}, there is a canonical embedding 
of the group of unipotent automorphisms of $R$ into 
the set of certain ``bifinite'' unipotent automorphisms
of the corresponding completion of $\calT$, and by the 
rigidity result \cite[Theorem 3.6]{Y2}, the latter are only coming from the 
center of $\calT$. This method puts very strong restrictions 
on the possible forms of the automorphisms of $R$ and with its  
help the problem of classifying $\Aut(R)$ can be treated 
with the currently developed ring theoretic techniques for 
studying $\Spec (R)$, see \cite{Y2,Y3}.

Embeddings of the form \eqref{eemb} are currently only known for 
quantum Schubert cell algebras and quantum double Bruhat cell algebras 
\cite{BZ,Ki,GLSh,GeY}. Those are derived using the Drinfeld 
$R$-matrix commutation relations for quantum function algebras. 
There is no general technique for constructing such embeddings 
for axiomatically defined families of algebras. 

For a very general 
(axiomatically defined) family of iterated skew polynomial extensions $R$, containing many quantum function algebras and quantized Weyl algebras, Cauchon \cite{Ca1} 
constructed embeddings of $R$ into quantum tori $\calT$
using the method of deleting derivations, which consists of formally 
exponentiating skew derivations in certain localizations of $R$. 
For those quantum tori the first embedding in \eqref{eemb}
is very rarely satisfied for the quantum affine space algebra $\calA$ 
generated by the corresponding Cauchon elements. The algebras which 
will play a key role in this paper were named  ``Cauchon--Goodearl--Letzter (CGL) extensions" in \cite{LLR}. Such an algebra is an iterated 
skew polynomial extension 
\begin{equation}
\label{Riter}
R = \KK [x_1] [x_2; \sig_2, \de_2] \cdots [x_N; \sig_N, \de_N],
\end{equation}
equipped with a rational action of 
a $\KK$-torus $\HH$ by algebra automorphisms,
which satisfies the following conditions:
\begin{enumerate}
\item[(i)] For all $1 \leq j < k \leq N$, $\sig_k(x_j) = \la_{kj} x_j$
for some $\la_{kj} \in \kx$.
\item[(ii)] For every $k \in [2,N]$, $\de_k$ is a locally nilpotent 
$\sig_k$-derivation of the subalgebra $R_{k-1}$
of $R$ generated by $x_1, \dots, x_{k-1}$.
\item[(iii)] The elements $x_1, \ldots, x_N$ are $\HH$-eigenvectors.
\item[(iv)] For every $k \in [1,N]$, there exists $h_k \in \HH$ such that $(h_k \cdot)|_{R_{k-1}} = \sig_k$ and
$h_k \cdot x_k = \la_k x_k$ for some $\la_k \in \kx$, which is 
not a root of unity.
\end{enumerate}

At this point we note that the subalgebras $\UU^\pm[w]$ of quantum Kac--Moody 
algebras $\UU_q(\g)$ introduced by De Concini--Kac--Procesi \cite{DKP} and Lusztig \cite{L} are iterated skew polynomial extensions which  satisfy the Levendorskii--Soibelman straightening rule and are CGL 
extensions with respect to an action of a torus arising from the 
root lattice grading of $\UU_q(\g)$. (See \exref{quantSchubert} for details.)

In this paper, we construct embeddings of the form \eqref{eemb} for all  CGL extensions $R$ using techniques from noncommutative 
unique factorization domains and a generalized version of the Gelfand--Tsetlin
procedure. In its original form, the Gelfand--Tsetlin procedure is used to construct 
large commutative subalgebras of noncommutative algebras
(usually universal enveloping algebras) or of Poisson algebras. 
One starts 
with a chain of $\KK$-algebra embeddings
\begin{equation}
\label{chain}
R_1 \subsetneq R_2 \subsetneq \cdots \subsetneq R_N = R.
\end{equation}
The Gelfand--Tsetlin subalgebra ${\mathcal{GT}}(R)$ of $R$ associated 
to this chain is the subalgebra of $R$ generated by 
$\ZZ(R_1) \cup \cdots \cup \ZZ(R_N)$.  It is obviously 
a commutative subalgebra of $R$. 

If the algebras $R_k$ are not universal enveloping algebras, 
their centers can be very small, yet they may 
have large monoids of normal elements. We consider 
the subalgebras generated by the latter -- for example,
in the case of quantum groups, those played a key role
in the classification of their maximal spectra \cite[Sections 5-7]{Y-sqg}.
For every $\KK$-algebra $R$, let us denote by $\NN(R)$ 
the subalgebra generated by the normal elements of $R$. We will 
call $\NN(R)$ the {\em{normal subalgebra}} of $R$.
With its help we introduce a twist to the 
Gelfand--Tsetlin construction. 
For each chain \eqref{chain}, we define the corresponding
{\em{normal Gelfand--Tsetlin subalgebra}} $\NGZ(R)$ of $R$   
as the subalgebra of $R$ generated by $\NN(R_1) \cup \cdots \cup \NN(R_N)$. 
(Note that $\NGZ(R)$ depends not only on $R$ but also 
on the choice of the chain \eqref{chain}. For simplicity, we will suppress this dependence in the 
notation.) Given an iterated skew polynomial extension \eqref{Riter}, 
there is a canonical choice for a chain of subalgebras \eqref{chain} 
given by the intermediate algebras in \eqref{Riter}.
In this setting, we prove the following result (\thref{embed}):

\bth{1} For all CGL extensions $R$, $\NGZ(R)$ is a quantum
affine space algebra, and \eqref{eemb} is satisfied for
$\calA : = \NGZ(R)$ and the quantum torus $\calT$ obtained by inverting 
the generators of $\calA$.
\eth

Noncommutative {\em{unique factorization domains}} were defined and studied 
by Chatters and Jordan \cite{Cha,ChJo}. They are domains $R$ with the property that 
each nonzero prime ideal of $R$ contains a prime element (a nonzero normal element of $R$ 
which generates a completely prime ideal). 
Assume that $R$ is equipped 
with a rational action of a torus $\HH$, or equivalently, $R$ is equipped 
with an $\xh$-grading, where $\xh$ is the character lattice of $\HH$.
The algebra $R$ is called an \emph{$\HH$-UFD} if each nonzero $\HH$-prime ideal of $R$ 
contains an $\xh$-homogeneous prime element. For a noetherian $\HH$-UFD $R$, 
$\NN(R)$ is precisely the subalgebra of $R$ generated by 
its $\xh$-homogenous prime elements, see Section \ref{NCUFD}
for details. By a theorem of Launois, Lenagan, and Rigal \cite{LLR}, 
all CGL extensions $R$ are noetherian $\HH$-UFDs. Hence, the 
normal Gelfand--Tsetlin subalgebra of a CGL extension $R$ is 
the subalgebra of $R$ generated by the set of $\xh$-homogeneous
prime elements of all intermediate algebras $R_1, \ldots, R_N$ from the iterated extension \eqref{Riter}.
We deduce \thref{1} from the following theorem which establishes 
a recursive relationship between the sets of $\xh$-homogeneous 
prime elements of the $\HH$-UFDs $R_1, \ldots, R_N$ for an arbitrary
CGL extension $R$. This is the major result in the paper (see \thref{CGL} and equation \eqref{y_comrel}).

For a natural number $n$, set 
$[1,n]:=\{1, \ldots, n\}$. Given a map $\mu : [1,N] \to [1,n]$, we 
define natural predecessor and successor functions 
$p : [1,N] \to [1,N] \sqcup \{ - \infty \}$ and
$s : [1,N] \to [1,N] \sqcup \{ + \infty \}$ for the level sets of $\mu$ by equations
\eqref{p} and \eqref{s}. 

\bth{2} For an arbitrary CGL extension $R$ of length $N$, there exist 
a natural number $n$, a surjective map $\mu : [1,N] \to [1,n]$ 
and a set of elements 
\[
\{ c_k \mid k \in [2,N], \; p(k) \ne -\infty \},
\]
such that $c_k \in R_{k-1}$ is $\xh$-homogeneous of the same degree 
as 
\[
\prod \{ x_j \mid j \in \mu^{-1} \mu(k) \cap [1, k] \}
\]
and the elements $y_1 \in R_1, \ldots, y_N \in R_N =R$ recursively given by
\[
y_k := 
\begin{cases}
y_{p(k)} x_k - c_k, &\mbox{if} \; \;  p(k) \neq - \infty \\
x_k, & \mbox{if} \; \; p(k) = - \infty
\end{cases}
\]
have the following properties:

{\rm(a)} For all $k \in [1,N]$, the $\xh$-homogeneous prime elements of $R_k$ are the elements 
\[
\kappa y_j, \; \; \kappa \in \kx, \; \; j \in [1, k], \; \;  s(j) > k.
\]

{\rm(b)} For all $1 \leq j < k \leq N$,
\[
y_k y_j = q_{kj} y_j y_k 
\]
where $q_{kj} \in \kx$ are certain explicit products of the elements
$\la_{kj} \in \kx$ given by equation \eqref{q}.
\eth

The number $n$ in \thref{2} equals  $| \{ k \in [1,N] \mid \de_k = 0 \} |$, by equation \eqref{rank}.

Using the method of \cite[Section 5]{Y-sqg}, we prove that each CGL extension $R$ is a free left and right $\NN(R)$-module in which $\NN(R)$
is a direct summand. We also construct an explicit $\NN(R)$-basis 
of $R$, see \thref{free} for details.

There is a torus going with any CGL extension $R$, which is maximal within a suitable subgroup of $\Aut(R)$, whose action covers all the automorphisms corresponding to normal elements of $R$, as follows (see Theorems \ref{tCGLfull}, \ref{tpinHmax} and \coref{normalHmax}). 

\bth{new3}
Let $R$ be an arbitrary CGL extension of length $N$, and let $n$ be the natural number appearing in Theorem {\rm\ref{t2}}. The group
\[
\HH' := \{ \psi \in \Aut(R) \mid x_1,\dots,x_N \text{\;are\;} \psi\text{-eigenvectors} \}
\]
is a $\KK$-torus of rank $n$, and $R$ is a CGL extension with respect to the action of $\HH'$. For every nonzero normal element $u \in R$, there exists $h \in \HH'$ such that $ua = (h\cdot a)u$ for all $a \in R$.
\eth

Along the way to the proof of \thref{2}, we establish a general 
recursive relationship between the set of $\xh$-homogeneous prime elements of 
an $\HH$-UFD $B$ and the sets of $\xh$-homogeneous prime elements of certain 
skew polynomial extensions $R:=B[x;\sig,\de]$. For this one 
needs to require that the conditions (Cx1)--(Cx5) listed 
in Section \ref{FromRtoA-a} are satisfied. These conditions are nothing but 
an abstraction of the conditions for one step of a CGL extension. 
Under those conditions, Launois, Lenagan, and Rigal proved 
\cite[Proposition 2.9]{LLR} that 
$R$ is also an $\HH$-UFD. The following result is \thref{RtoBa}.
 
\bth{3}
Let $R=B[x; \sig, \de]$ be a skew polynomial algebra, equipped with a rational action of a $\KK$-torus $\HH$ by algebra automorphisms which leave $B$ invariant. Assume that $B$ is an $\HH$-UFD
and that $R$ satisfies the conditions {\rm(Cx1)--(Cx5)} listed in Section {\rm\ref{FromRtoA-a}}.
Let $\{u_i \mid i \in I\}$ be a list of the $\xh$-homogeneous prime elements
of $B$ up to taking associates. Then there are the following three possibilities 
for a list of the $\xh$-homogeneous prime elements of $R$ up to taking associates:

{\rm (i)} $\{u_i \mid i \in I, i \neq i_0  \} \sqcup \{ v_{i_0}:=u_{i_0} x - c_\circ \}$
for some $i_0 \in I$ and $c_\circ \in B$ such that $u^{-1}_{i_0} c_\circ$ is a nonzero 
$\xh$-homogeneous element of $B[u^{-1}_{i_0}]$ with the same $\xh$-degree as $x$.

{\rm (ii)} $\{ u_i \mid i \in I \} \sqcup \{ x \}$.

{\rm (iii)} $\{ u_i \mid i \in I \}$.
\eth 

Section \ref{NCUFD} contains definitions and some general facts about noncommutative 
UFDs and prime elements. \thref{3} is proved in Section \ref{FromRtoA-a}. 
In Section \ref{CGL}, we prove that the situation (iii) never occurs in the 
setting of CGL extensions. This establishes a general vanishing 
property of certain skew derivations of CGL extensions, see \thref{CGLind}(b) 
for details. Theorems \ref{t1} and \ref{t2} are also proved in Section \ref{CGL}, while \thref{new3} is proved in Section \ref{fullness}.

For the iterated skew polynomial extensions $R$ which are subalgebras of quantized universal 
enveloping algebras $\UU_q(\g)$, the elements $x_1, \ldots, x_N$ are among the 
Lusztig root vectors of $\UU_q(\g)$. They satisfy the Levendorskii--Soibelman
straightening rule, which means that in the setting of \eqref{Riter}, 
$\de_k(x_j)$ belongs to the subalgebra of $R$ generated by 
$x_{j+1}, \ldots, x_{k-1}$ for all $j < k$. Together with a mild assumption on the action of $\HH$,
all such CGL extensions $R$ have a second CGL extension presentation
\begin{equation}
\label{second}
R = \KK[x_N] [x_{N-1}; \sig'_{N-1},\de'_{N-1}] \cdots [x_1; \sig'_1, \de'_1].
\end{equation}
We call such CGL extensions \emph{symmetric}, see \deref{symmetric}.
In \thref{rel}, we prove that for all symmetric CGL extensions
$R$, the quantum torus associated to $\NGZ(R)$ for the presentation \eqref{Riter} 
coincides with the Cauchon quantum torus \cite{Ca1} of $R$ for the presentation
\eqref{second}. We also explicitly express the $\xh$-homogeneous prime elements of all 
intermediate algebras $R_1, \ldots, R_N$ in terms of the Cauchon 
elements of $\Fract(R)$ (for the presentation \eqref{second}). 
The proof is based on an extension of the method of \cite[Section 3]{GeY}.  

Commutative unique factorization domains were previously used in the
area of cluster algebras \cite{BFZ,GLSh2}. The idea was that after enough
clusters are constructed for certain coordinate rings, the unique
factorization property can be used to prove that the coordinate rings
actually coincide with the constructed (upper) cluster algebras. In
our treatment, the noncommutative UFD property is used in a different
fashion, namely to construct initial clusters for a more general
family of noncommutative associative algebras which do not necessarily
come from Kac--Moody Lie algebras.
We prove in \cite{GYmem} that, under mild assumptions on the scalars in the base field, every symmetric CGL extension possesses a quantum cluster algebra structure. (See \cite{GYann} for an overview.)

\subsection{Notation}
\label{notation}
We finish the introduction with a word on notation. Throughout the paper, $\KK$ 
will denote an infinite base field. All algebras will be unital $\KK$-algebras, 
all automorphisms will be $\KK$-algebra automorphisms, and all skew derivations 
will be $\KK$-linear left skew derivations. All skew polynomial rings $B[x;\sig,\de]$ will be assumed to be left Ore extensions, meaning that the commutation rule is given by $xb = \sig(b) x + \de(b)$ for $b\in B$. Elements $a$ and $b$ in a $\KK$-algebra are said to \emph{quasi-commute} if $ab$ and $ba$ are nonzero scalar multiples of each other. The center of a ring $R$ is denoted $\ZZ(R)$.

\sectionnew{ Noncommutative unique factorization domains}
\label{NCUFD}

We recall the noncommutative unique factorization conditions introduced by Chatters and Jordan \cite{Cha}, \cite{ChJo}, summarize some of their results and a theorem of Launois, Lenagan, and Rigal \cite{LLR}, and develop some extensions.

\subsection{UFRs and UFDs}
\label{UFRD}
We separate out the noetherian assumption from \cite[Definition, p.50]{Cha}, \cite[Definition, p.23]{ChJo}, and \cite[Definitions 1.1, 1.2]{LLR}.

A (noncommutative) {\em{unique factorization ring}} ({\em{UFR}}) is a prime ring $R$ such that each nonzero prime ideal of $R$ contains a nonzero prime ideal generated by a normal element, i.e., an element $u\in R$ such that $Ru= uR$. 

Let $R$ be a domain. A {\em{prime element}} in $R$ is any nonzero normal element $p\in R$ such that $Rp$ is a completely prime ideal, i.e., $R/Rp$ is a domain. A (noncommutative) {\em{unique factorization domain}} ({\em{UFD}}) is a domain $R$ such that each nonzero prime ideal of $R$ contains a prime element.

Because of the noncommutative Principal Ideal Theorem (e.g., \cite[Theorem 4.1.11]{McRob}), a noetherian UFR is a UFD if and only if it is a domain and all its height one prime ideals are completely prime.

Continue to assume that $R$ is a domain, and let $u\in R$ be a normal element. We say that $u$ is a {\em{divisor}} of an element $a\in R$, written $u\mid a$, provided $a\in Ru$, which holds if and only if $a=ru$ for some $r\in R$, if and only if $a=us$ for some $s\in R$. As in the commutative case, a nonzero normal element $p\in R$ is prime if and only if $p$ is not a unit and $(p\mid ab \implies p\mid a$ or $p\mid b)$ for all $a,b\in R$. We say that normal elements $u,v \in R$ are {\em{associates}} provided $Ru= Rv$, which occurs if and only if $u= av$ for some unit $a\in R$, if and only if $u=vb$ for some unit $b\in R$.

\bpr{factorUFD}
Let $R$ be a noetherian UFR.

{\rm(a)} Every nonzero normal element of $R$ can be expressed in the form $cp_1p_2 \cdots p_n$ for some unit $c\in R$ and some normal elements $p_i\in R$ such that each $Rp_i$ is a prime ideal. In case $R$ is a UFD, the $p_i$ must be prime elements.

Now assume that $R$ is a noetherian UFD.

{\rm(b)} Every nonzero element of $R$ can be expressed in the form $cp_1p_2 \cdots p_n$ for some prime elements $p_i\in R$ and some $c\in R$ which has no prime divisors.

{\rm(c)} Let $u$ be a nonzero, nonunit, normal element of $R$. Then $u$ is prime if and only if $u$ is \underbar{irreducible}, that is, {\rm(}$u=ab \implies a$ or $b$ is a unit{\rm\/)} for all $a,b\in R$.

{\rm(d)} Let $u\in R$ be a nonzero normal element, and let $a,b\in R$ such that $u\mid ab$ or $u\mid ba$. If $u$ and $a$ have no prime common divisors, then $u\mid b$.
\epr

\begin{proof} For (a) and (b), see \cite[p.24]{ChJo} and \cite[Proposition 2.1]{Cha}. Parts (c) and (d) follow from (b) just as in the commutative case.
\end{proof}

\subsection{$\HH$-UFDs}
\label{HUFD}
Let $R$ be a domain which is a $\KK$-algebra, and $\HH$ a group acting on $R$ by $\KK$-algebra automorphisms. Recall that an {\em{$\HH$-prime ideal}} of $R$ is any proper $\HH$-stable ideal $P$ of $R$ such that $(IJ\subseteq P \implies I\subseteq P$ or $J\subseteq P)$ for all $\HH$-stable ideals $I$ and $J$ of $R$.

Following \cite[Definition 2.7]{LLR}, we say that $R$ is an {\em{$\HH$-UFD}} provided each nonzero $\HH$-prime ideal of $R$ contains a prime $\HH$-eigenvector. The proof of \prref{factorUFD} 
(see also \cite[Proposition 6.18 (ii)]{Y-sqg}) 
is easily adapted to give the following result.

\bpr{factorHUFD}
Let $R$ be a noetherian $\HH$-UFD.

{\rm(a)} Every normal $\HH$-eigenvector in $R$ is either a unit or a product of prime $\HH$-eigen\-vectors.

{\rm(b)} Let $u$ be a nonunit normal $\HH$-eigenvector in $R$. Then $u$ is prime if and only if it is irreducible. 

{\rm(c)} Let $u\in R$ be a normal $\HH$-eigenvector, and let $a,b\in R$ such that $u\mid ab$ or $u\mid ba$. If $u$ and $a$ have no prime $\HH$-eigenvector common divisors, then $u\mid b$.
\epr

\bre{notation}  Let $R$ 
be a $\KK$-algebra with a rational action of a $\KK$-torus $\HH$ 
by $\KK$-algebra automorphisms. Equivalently, $R$ is equipped with a grading by the character lattice $\xh$ of $\HH$ (e.g., \cite[Lemma II.2.11]{BG}). An element of $R$ is an $\HH$-eigenvector if and only if
it is nonzero and $\xh$-homogeneous. Throughout the paper we will use the term 
{\em{homogeneous prime element}} of $R$ instead of {\em{prime $\HH$-eigenvector}}
since this is more suggestive from a ring theoretic perspective.
The term {\em{homogeneous}} will always mean $\xh$-{\em{homogeneous}}.
The latter term will only be used in cases when the dependence on
the underlying torus $\HH$ has to be emphasized. 

Note that $\xh$ is a free abelian group and therefore a totally ordered group. Consequently, if $R$ is a domain, then all units of $R$ are homogeneous. Finally, we recall from \cite[Proposition II.2.9]{BG} that if $R$ is noetherian, all $\HH$-prime ideals of $R$ are prime. 
\ere

The arguments of \cite{LLR} yield the following result.

\bth{HUFD-UFR}
{\rm [Launois-Lenagan-Rigal]}
Let $R$ be a noetherian $\KK$-algebra, and $\HH$ a $\KK$-torus acting rationally on $R$ by $\KK$-algebra automorphisms. Assume that $R$ is an $\HH$-UFD. Then $R$ is a UFR, but not necessarily a UFD.
\eth

\begin{proof} (Adapted from \cite[Proposition 1.6, Theorem 3.6]{LLR}.)

Let $X_0$ be the set of homogeneous prime elements in $R$, and let $X$ be the multiplicative set generated by $X_0$. Since $X$ consists of normal elements, it is a right and left denominator set in $R$. Let $T= RX^{-1}$ be the corresponding Ore localization of $R$.
We first show that any nonzero $\HH$-stable ideal $I$ of $R$ meets $X$. Since $R$ is noetherian, there are prime ideals $P_1,\dots,P_t$ minimal over $I$ such that $P_1P_2\cdots P_t \subseteq I$. For each $j$, the largest $\HH$-stable ideal $Q_j$ contained in $P_j$ is $\HH$-prime, and it is nonzero because it contains $I$. By hypothesis, $Q_j$ contains a homogeneous prime element $q_j$. Therefore $I$ contains the element $q_1q_2 \cdots q_t$ from $X$.

It follows that $T$, with respect to the induced action of $\HH$, is an {\em{$\HH$-simple}} ring, i.e., the only $\HH$-stable ideals of $T$ are $0$ and $T$. By \cite[Corollary II.3.9]{BG}, the center of $T$ is a Laurent polynomial ring over a field and there are inverse bijections between $\Spec T$ 
and $\Spec \ZZ(T)$ given by contraction and extension. Hence, $\ZZ(T)$ is a commutative UFD, and every nonzero prime ideal of $T$ contains a nonzero prime ideal generated by a central element.

We must show that each nonzero prime ideal $P$ of $R$ contains a nonzero prime ideal generated by a normal element. If $P$ meets $X$, then because the elements of $X_0$ are normal, $P$ must meet $X_0$. Take $p\in P\cap X_0$; then $P$ contains the nonzero prime ideal $Rp$.

Now assume that $P\cap X= \varnothing$. Then $PT$ is a nonzero prime ideal of $T$, so there is some nonzero element $z\in PT\cap \ZZ(T)$ such that $zT$ is a prime ideal. The contraction $P'= zT\cap R$ is a nonzero prime ideal of $R$ such that $P'\subseteq P$ and $P'T= zT$, and we may replace $P$ by $P'$. Thus, there is no loss of generality in assuming that $PT= zT$.  We can write $z= ax^{-1}$ for some $a\in P$ and $x\in X$, and we have $P= a_1R+ \cdots+ a_mR$ for some $a_i\in P$. Moreover, each $a_i= zt_i$ for some $t_i\in T$. There exist $p_1,\dots,p_n \in X_0$ such that $x^{-1}$, $z$, $t_1,\dots,t_m$ all lie in $RY^{-1}$, where $Y$ is the multiplicative set generated by $\{p_1,\dots,p_n\}$. Since we may remove any $p_j$ which is an associate of a different $p_i$, we may assume that $p_i$ is not an associate of $p_j$ for any $i\ne j$. Note that $z\in \ZZ(RY^{-1})$ and $P(RY^{-1})= zRY^{-1}$.

Set $R_0=R$ and $R_i= R[p_1^{-1},\dots,p_i^{-1}]$ for $i=1,\dots,n$. As in the proof of \cite[Proposition 1.6]{LLR}, it follows from \cite[Lemma 1.5]{LLR} that $p_{i+1}$ is a prime element of $R_i$ for each $i=0,\dots,n-1$. A reverse induction using \cite[Lemma 1.4]{LLR} now shows that $PR_i$ is generated by a normal element of $R_i$ for each $i=n,\dots,0$. Therefore $P=Ru$ for some normal element $u$ in $R$. This completes the proof that $R$ is a UFR.

To see that $R$ need not be a UFD, assume that $\chr\KK \ne 2$ and  take $R$ to be the quantum plane $\KK\langle x,y\mid xy=-yx\rangle$. The torus $\HH= (\kx)^2$ acts rationally on $R$ so that $(\alpha,\beta).x= \alpha x$ and $(\alpha,\beta).y = \beta y$ for all $(\alpha,\beta) \in \HH$. The elements $x,y\in R$ are homogeneous and prime, and any nonzero $\HH$-prime ideal of $R$ contains one of them (e.g., \cite[Example II.1.11]{BG}). Thus, $R$ is an $\HH$-UFD. The element $x^2-1$ is central in $R$, and it is easily seen that $R(x^2-1)$ is a prime ideal (e.g., apply \cite[Lemma 10.6.4(iv)]{McRob}), necessarily of height one. It is not completely prime, and therefore $R$ is not a UFD.
\end{proof}

Now return to the situation at the beginning of this subsection.
A normal element $u\in R$ is said to be {\em{$\HH$-normal}} if there is some $h\in \HH$ such that $ua= h(a)u$ for all $a\in R$. (It is not assumed that $u$ is an $\HH$-eigenvector.) We reserve 
the term {\em{$\HH$-prime element}} for any prime element of $R$ which is both $\HH$-normal and an $\HH$-eigenvector, i.e., any $\HH$-normal homogeneous prime element. The algebra $R$ will be called a {\em{strong $\HH$-UFD}} in case each nonzero $\HH$-prime ideal of $R$ contains an $\HH$-prime element. An $\HH$-UFD $R$ is a strong $\HH$-UFD if and only if each homogeneous prime element of $R$ 
has an associate which is an $\HH$-prime element. If $R$ is an $\HH$-UFD and all 
units of $R$ are central, then $R$ is a strong $\HH$-UFD if and only if its homogeneous prime elements are $\HH$-prime.

\bpr{factorstrongHUFD}
Let $R$ be a noetherian strong $\HH$-UFD.

{\rm(a)} Every normal $\HH$-eigenvector in $R$ can be expressed in the form $cp_1p_2\cdots p_n$ for some unit $\HH$-eigenvector $c\in R$ and some $\HH$-prime elements $p_i\in R$.

{\rm(b)} Let $u\in R$ be a normal $\HH$-eigenvector, and let $a,b\in R$ such that $u\mid ab$ or $u\mid ba$. If $u$ and $a$ have no $\HH$-prime common divisors, then $u\mid b$.

{\rm(c)} Let $a$ and $b$ be nonzero elements of $R$. Then $a=a'w$ and $b=b'w$ where $w$ is a product of $\HH$-prime elements of $R$ while $a'$ and $b'$ are elements of $R$ with no $\HH$-prime common divisors.
\epr

\begin{proof} Part (a) is proved in the same way as in Propositions \ref{pfactorUFD}, \ref{pfactorHUFD}, and parts (b), (c) follow from (a) just as in the commutative case.
\end{proof}

\bth{HnormalsinHUFD}
Let $R$ be a noetherian $\KK$-algebra, and $\HH$ a $\KK$-torus acting rationally on $R$ by $\KK$-algebra automorphisms. Assume that $R$ is a strong $\HH$-UFD. Then each normal element of $R$ is an associate of an $\HH$-normal element. If all units of $R$ are central, then all normal elements of $R$ are $\HH$-normal.
\eth

\begin{proof} The second conclusion of the theorem is an immediate consequence of the first. By \thref{HUFD-UFR}, $R$ is a UFR. Thus, to prove the first conclusion of the present theorem, \prref{factorUFD}(a) allows us to reduce to the case of a normal element $p\in R$ such that $Rp$ is a prime ideal. If $p$ is an associate of an $\HH$-eigenvector, then $Rp$ is a height one $\HH$-prime ideal of $R$. In this case, $Rp=Rp'$ for some $\HH$-prime element $p'$, and we are done. Thus, we may assume that $p$ is not an associate of any $\HH$-prime element. Since $Rp$ is a height one prime ideal of $R$, it thus cannot contain any $\HH$-prime elements.

Let $X$ be the multiplicative set generated by the set of $\HH$-prime elements of $R$, and set $T= RX^{-1}$. As in the proof of \thref{HUFD-UFR}, $\ZZ(T)$ is a commutative UFD, and contraction and extension provide inverse bijections between $\Spec T$ and $\Spec \ZZ(T)$. Since $Rp$ is a prime ideal containing no $\HH$-prime elements, and these elements are all normal, it follows that $Rp$ is disjoint from $X$. Consequently, $Tp$ is a height one prime ideal of $T$. Thus, $Tp= Tp'$ for some $p'\in \ZZ(T)$, and $p=p't$ for some unit $t\in T$. We claim that $t= cu_1u_2^{-1}$ for some unit $c$ of $R$ and some $u_1,u_2\in X$. Then $c^{-1}pu_2= p' u_1$, and there exist $h_1,h_2\in \HH$ such that $u_ia= h_i(a)u_i$ for all $a\in R$. Now
\[
c^{-1}ph_2(a)u_2= c^{-1}pu_2a= p' u_1a= h_1(a)p' u_1= h_1(a)c^{-1}pu_2
\]
for all $a\in R$, and thus $c^{-1}pb= h_1h_2^{-1}(b)c^{-1}p$ for all $b\in R$, proving that $c^{-1}p$ is $\HH$-normal.

It remains to show that any unit $t\in T$ must have the claimed form. Write $t= au^{-1}$ and $t^{-1}= v^{-1}b$ for some $a,b\in R$ and $u,v\in X$. In view of \prref{factorstrongHUFD}(c), we may assume that $a$ and $u$ have no $\HH$-prime common divisors, and that $b$ and $v$ have no $\HH$-prime common divisors. From the equation $v^{-1}bau^{-1}=1$ in $T$, we obtain $ba= vu$ in $R$. By \prref{factorstrongHUFD}(b), $u\mid b$ and $v\mid a$, say $b=ub'$ and $a=a'v$ for some $a',b'\in R$. Since $u$ and $v$ quasi-commute (because both are $\HH$-normal and homogeneous), we conclude that $b'a'$ is a nonzero scalar. Consequently, $a'$ is a unit in $R$, and therefore $t= a'vu^{-1}$ has the desired form.
\end{proof}

Recall the definition of the normal subalgebra $\NN(R)$ of a 
$\KK$-algebra $R$. We have the following fact concerning the normal 
subalgebra of an $\HH$-UFD $R$.

\bpr{norm-UFD} Let $R$ be a noetherian $\HH$-UFD for a $\KK$-torus
$\HH$ acting rationally on $R$ by $\KK$-algebra automorphisms. Then $\NN(R)$ 
is precisely the subalgebra of $R$ generated by the 
homogeneous prime elements and the units of $R$.
\epr

If the group of units of $R$ is reduced to scalars, 
then \prref{norm-UFD} states that $\NN(R)$ is the (unital) subalgebra 
of $R$ generated by the homogeneous prime elements of $R$.

\begin{proof} Obviously the subalgebra of $R$ generated 
by the homogeneous prime elements and units of $R$ is a subalgebra of $\NN(R)$. 
For the opposite inclusion, let $u \in R$ be a nonzero 
normal element. 
By \cite[Proposition 6.20]{Y-sqg},
\[
u = u_1 + \cdots + u_n
\]
where $u_1, \ldots, u_n$ are nonzero normal $\HH$-eigenvectors
with distinct eigenvalues
and by \prref{factorHUFD}(a) each of them is either a unit 
or a product of homogeneous prime elements in $R$.
\end{proof}

\sectionnew{From prime elements of $B$ to prime elements of $B[x;\sig, \de]$}
\label{FromRtoA-a}

\subsection{General assumptions}
\label{2.1a}
Throughout Section \ref{FromRtoA-a}, we work in the following setting:
\begin{enumerate}
\item[$\bullet$] {\em{$B$ is a $\KK$-algebra with a rational action of a $\KK$-torus $\HH$
by $\KK$-algebra automorphisms.}}
\item[$\bullet$] {\em{$B$ is an $\HH$-UFD {\rm(}recall Subsection {\rm\ref{HUFD})}.}}
\item[$\bullet$] {\em{$R=B[x; \sig, \de]$ is a skew polynomial algebra, equipped with a rational action of $\HH$ by $\KK$-algebra automorphisms, extending the action of $\HH$ on $B$.}}
\end{enumerate}

Consider the following conditions:
\begin{enumerate}
\item[(Cx1)] {\em{$B$ is noetherian}}.
\item[(Cx2)] {\em{$\delta$ is locally nilpotent}}.
\item[(Cx3)] {\em{$B$ is $\HH$-stable and $x$ is an $\HH$-eigenvector. If this holds, let $\la\in\xh$ denote the $\HH$-eigenvalue of $x$}}.
\item[(Cx4)] {\em{There exists $h_\circ \in \HH$ such that $(h_\circ \cdot)|_B = \sig$ and $x$ is an $h_\circ$-eigenvector with $h_\circ$-eigen\-value $\la_\circ$ which is not a root of unity}}.
\item[(Cx5)] {\em{All $\HH$-prime ideals of $B$ are completely prime}}.
\end{enumerate}
Of course, if (Cx1) holds, then $R$ is noetherian as well.

Conditions (Cx1)--(Cx5), together with the assumption that $B$ is a $\KK$-algebra domain, mean that  $R$ is a {\em{Cauchon extension}} in the sense of \cite[Definition 2.5]{LLR}. (The $q$-skew condition in that definition holds automatically, due to the following observation.) If (Cx4) holds, then, by applying $h_\circ$ to the equations $xa=\sig(a)x +\de(a)$ for $a\in B$, we see that $\la_\circ x\sig(a)= \sig^2(a)\la_\circ x+ \sig\de(a)$; comparing this with $x\sig(a)= \sig^2(a) x+ \de\sig(a)$, we conclude that 
\begin{equation} \label{sigdel}
\sig\de(b)= \la_\circ \de\sig(b) \qquad \forall b\in B.
\end{equation}
Similarly, if (Cx3) holds, we see that
\begin{equation} \label{hdel}
(h\cdot)|_B \de= \la(h)\de(h\cdot)|_B \qquad \forall h\in\HH.
\end{equation}

By \cite[Proposition 2.9]{LLR}, if conditions 
{\rm{(Cx1)--(Cx5)}} are satisfied, then $R$ is an $\HH$-UFD. 
The goal of this section is to obtain an explicit description 
of the set of homogeneous prime elements of $R$ in terms of the 
homogeneous prime elements of $B$.

If $u \in B$ is a nonzero normal element, then 
the corresponding automorphism of $B$ will be denoted by 
\begin{equation}
\label{aut-u}
\varphi_u \in \Aut (B), \; \; 
\mbox{where} \; \; 
u b = \varphi_u(b) u, \qquad \forall b \in B.
\end{equation}

\subsection{Degree one homogeneous prime elements of $R$ and properties of the skew derivation $\de$}
\label{2.2a}
Our minimal assumptions in this subsection are conditions (Cx3) and (Cx4).
Denote by $E$ the multiplicative set generated by all homogeneous prime elements of $B$. The set $E$ is an Ore subset of $B$. Since $\sig(E)= h_\ci \cdot E  = E$, this set is also an Ore subset of $R$, and 
\[
R[E^{-1}] = B[E^{-1}][x;\sig,\de]
\]
for the canonical induced actions of $\sig$ and $\de$ on $B[E^{-1}]$.
Since $B$ is an $\HH$-UFD, every nonzero $\HH$-prime ideal of $B$
meets $E$. In case $B$ is noetherian, it follows (as in the proof of \thref{HUFD-UFR}) that every nonzero $\HH$-stable ideal of $B$ meets $E$, and thus $B[E^{-1}]$ is $\HH$-simple. 

We will need the following result which is a special case of 
\cite[Lemma II.5.10]{BG}:

\ble{nonHsimp} Assume {\rm(Cx1), (Cx3)} and {\rm(Cx4)}. 
If in the above setting $R[E^{-1}]$ is not $\HH$-simple, then there 
exists a unique homogeneous element $d\in B[E^{-1}]$ such that $d=0$ or $d$ 
has the same $\xh$-degree as $x$ and 
\[
\de(b)= d b - \sig(b)d \; \;  
\mbox{for all} \; \;  
b \in B[E^{-1}].
\] 
Moreover, $R[E^{-1}](x-d)$ is the only nonzero 
$\HH$-prime ideal of $R[E^{-1}]$.
\ele

\bco{nondeg0} Assume that {\rm(Cx1), (Cx3)} and {\rm(Cx4)}
are satisfied. Then all homogeneous prime elements of $R = B[x;\sig,\de]$ 
have degree at most {\rm1} in $x$. Up to taking associates, 
there is at most one homogeneous prime element of $R$ 
which does not belong to $B$ {\rm(}i.e., has degree {\rm1} in $x${\rm)}.
\eco

\begin{proof} If $v$ is a homogeneous prime element of $R$ such that 
$v \notin B$, then $v$ is a homogeneous prime element of $R[E^{-1}]$. 
Thus, $R[E^{-1}]$ is not $\HH$-simple, so there exists $d\in B[E^{-1}]$ as in \leref{nonHsimp}, and $v$ is an 
associate of the prime element $x - d \in R[E^{-1}]$ (as prime 
elements of $R[E^{-1}]$). This implies that $v$ has degree at most $1$ in $x$. Since $Rv$ is a prime ideal of $R$ disjoint from $E$, we have
\[
Rv = (R[E^{-1}]v) \cap R = 
\big( R[E^{-1}] (x-d) \big) \cap R.
\]

If $w$ is any other homogeneous prime element of $R$ that is not in $B$, the same argument as above shows that
\[
Rw = \big( R[E^{-1}] (x-d) \big) \cap R,
\]
and therefore $Rw= Rv$.
\end{proof}

We finish this subsection with some properties of locally nilpotent 
skew derivations of an $\HH$-UFD.

\ble{innerdelta-a} Assume {\rm(Cx2)--(Cx4)}.

{\rm (a)} If $a\in B$ such that $\de(a)\in aB$ or $\de(a) \in Ba$, then $\de(a) = 0$.

{\rm (b)} If $\de$ is an inner $\sig$-derivation, i.e., 
there is some $c \in B$ such that $\de(b) = c b - \sig(b) c$, $\forall b \in B$, 
then $\de=0$. If, in addition $c$ is a 
homogeneous element of $B$ of degree $\la$, then $c=0$.
\ele

\begin{proof} (a) By assumption, $\de$ is locally nilpotent and $\sig\de = \la_\circ \de\sig$ with $\la_\circ$ not a root of unity. The result is thus given by \cite[Lemme 7.2.3.2]{Ric}.

(b) Write $c= c_0+ \cdots+ c_m$ where the $c_i$ are homogeneous elements (possibly zero) of $B$ 
of distinct $\xh$-degrees $\mu_i$ and $\mu_0= \la$. 
If $b \in B$ is homogeneous of degree $\rho$, then
$\de(b)$ is homogeneous of degree $\la\rho$ by \eqref{hdel}.
On the other hand, each $c_i b- \sig(b)c_i$ is homogeneous of degree $\mu_i\rho$, 
so $\de(b)= c_0 b - \sig(b)c_0$. Since this holds for all homogeneous 
$b \in B$, it holds for all $b \in B$.
Now $\de(c_0)= c_0^2- (h_\circ \cdot c_0)c_0= (1-\la_\circ) c_0^2$. Part (a) implies that $\de(c_0) = 0$. Since $\la_\circ \ne 1$, it follows that $c_0^2 = 0$ and so
$c_0=0$. Thus, $\de=0$.
\end{proof}

\bre{moregen} Although we will not need this, we note that 
the second statement in \leref{innerdelta-a}(b) holds under the weaker 
assumption that $c$ is a $\sig$-eigenvector. The proof 
of this fact is the same.
\ere

\subsection{Inducing prime elements of $R=B[x;\sig,\de]$ from prime elements of $B$}
\label{2.4a} 
Next we construct a homogeneous prime element of $R=B[x;\sig,\de]$ for each 
homogeneous prime element of $B$.

Assume that $R=B[x;\sig, \de]$ and $\HH$ satisfy the conditions 
{\rm{(Cx1)--(Cx4)}}. Then the set $\{ x^n \mid n \in \Znn \}$ is an 
Ore subset of $R$ (e.g., \cite[Lemme 2.1]{Ca1}). Consider the Cauchon map
\begin{equation}
\label{C}
\theta : B \to R[x^{-1}] \; \; 
\mbox{given by} \; \; 
\theta(b)= 
\sum_{m=0}^\infty \frac{(1- \la_\circ)^{-m}}{[m]_{\la_\circ}!} [ \delta^m \sig^{-m}(b)] 
x^{-m},
\end{equation}
where $\la_\circ \in \kx$ is the element from {\rm{(Cx4)}} and $[m]_q = 1+ \cdots + q^{m-1}, 
[m]_q!= [0]_q \cdots [m]_q$ 
are the standard $q$-integers and factorials. The $\HH$-action on $R$ induces an 
$\HH$-action on $R[x^{-1}]$ and the map $\theta$ is $\HH$-equivariant \cite[Lemma 2.6]{LLR}.
Denote $B':=\theta(B)$. 
Cauchon \cite[Section 2]{Ca1} proved that $\theta : B \to B'$ is a $\KK$-algebra isomorphism
and that $R[x^{-1}] = B' [x^{\pm 1}; \sig]$ for the extension of 
$\sig$ to an automorphism of $R[x^{-1}]$ (preserving $B'$) 
given by $\sig = (h_\circ \cdot)$.

For an ideal $J$ of $R$ denote the ideal of its leading coefficients
\[
\lc(J) := \{ b \in B \mid \exists  a \in J, \; m \in \Znn \; \; 
\mbox{such that} \; \; a -  b x^m \in B x^{m-1} + \cdots + B \}.
\]

\ble{new-prime-a} Assume {\rm{(Cx1)--(Cx5)}}.
If $J$ is an $\HH$-invariant ideal of 
$B$, then $\theta(J)[x^{\pm 1}] \cap R$ is an $\HH$-invariant ideal of 
$R$ and $\lc( \theta(J)[x^{\pm 1}] \cap R)=J$. If $J$ is completely prime, 
then $\theta(J)[x^{\pm 1}] \cap R$ is completely prime. If 
$J$ is a height one prime ideal of $B$, then $\theta(J)[x^{\pm 1}] \cap R$ 
is a height one prime ideal of $R$.
\ele

\begin{proof} The proofs of the first two statements are straightforward 
and are left to the reader. 

For the proof of the last statement, assume that $J$ is a height one prime 
ideal of $B$ for which $\theta(J)[x^{\pm 1}] \cap R$ is not a  
height one prime ideal of $R$. Since $R$ is an $\HH$-UFD by 
\cite[Proposition 2.9]{LLR}, there exists a homogeneous prime element 
$p$ of $R$ such that $R p \subsetneq \theta(J)[x^{\pm 1}] \cap R$. 
Thus $P:= R p$ is an $\HH$-invariant completely prime ideal of $R$ such that 
\[
\{ 0 \} \subsetneq P \subsetneq \theta(J)[x^{\pm 1}] \cap R.
\]
Furthermore, $x \notin \theta(J)[x^{\pm 1}] \cap R$ since 
$1 \notin \lc( \theta(J)[x^{\pm 1}] \cap R) = J$, where we 
used the first part of the lemma. Localization with respect to $\{x^m \mid m \in \Znn \}$
leads to  
\begin{equation}
\label{em3}
\{ 0 \} \subsetneq P[x^{- 1}] \subsetneq ( \theta(J)[x^{\pm 1}] \cap R) [x^{-1}] = 
\theta(J)[x^{\pm 1}].
\end{equation}
But $R[x^{-1}] = B'[ x^{\pm 1} ; \sig]$. 
Since $P[x^{-1}]$ is
an $\HH$-invariant ideal of $B'[x^{\pm 1}; \sig]$, by \cite[Lemma 2.2]{LLR}
\[
P[x^{-1}] = (P[x^{-1}] \cap B' ) [x^{\pm 1}]
\]
and \eqref{em3} implies 
\[
\{ 0 \} \subsetneq P[x^{-1}] \cap B' \subsetneq \theta(J).
\]
This is a contradiction, since $P[x^{-1}] \cap B'$ is a completely prime ideal of $B'$, 
and $\theta(J)$ is a height one prime 
ideal of $B'$.
\end{proof}

\bth{R-B} Let $R=B[x; \sig, \de]$ be a skew polynomial algebra, equipped with a rational action of a $\KK$-torus $\HH$ by algebra automorphisms which leave $B$ invariant. Assume that $B$ is an $\HH$-UFD
and that $R$ satisfies {\rm(Cx1)--(Cx5)}. Let $u$ be a homogeneous prime element  of $B$, and let $\al \in \kx$ be such that 
\[
\sig(u) = h_\circ \cdot u = \al u.
\] 
Then exactly one of the 
following two situations occurs:

{\rm{(i)}} The element $u$ remains a prime element of $R$.

In this case, $u$ quasi-commutes with 
$x$ via
\[
u x = \al^{-1} x u.
\]
Furthermore,
\[
\de(u) =0 \; \; \mbox{and} \; \; 
\theta(u B) [x^{\pm 1} ] \cap R = u R.
\]

{\rm{(ii)}} There exists a unique element $c_\circ \in B$ such that 
$v:=ux - c_\circ$ is a homogeneous 
prime element of $R$ {\rm{(}}in particular $u \nmid c_\circ${\rm{)}}. 

In this case,
$\de$ is given by 
\begin{equation}
\label{pr-b}
\de(b) = (u^{-1} c_\circ) b - \sig(b ) (u^{-1} c_\circ), \; \; \forall b \in B,
\end{equation}
and $v$ normalizes the elements of $R$ as follows:
\begin{equation}
\label{comm-ii}
v x = \al^{-1} x v, \quad
v b =  (\varphi_u \sig(b)) v = \varphi_u( h_\circ \cdot b) v , \; \; 
\forall b \in B,
\end{equation}
cf. \eqref{aut-u}. Furthermore,
\begin{align*}
u c_\circ &= (\al \la_\circ)^{-1} c_\circ u,  &\de(u) &= \al(\la_\circ-1) c_\circ,  \\
\de(c_\circ) &= 0,  &\de( u^{-1} c_\circ) &= -(\la_\circ-1) (u^{-1} c_\circ)^2,
\end{align*}
and
\[
\theta(u B) [x^{ \pm 1}] \cap R = v R.
\]
\eth

We note that by \coref{nondeg0}, the situation (ii) cannot simultaneously 
occur for two homogeneous prime elements $v$ of $R$ which are not associates 
of each other.

\begin{proof} By \leref{new-prime-a}, 
$\theta( u B) [x^{\pm 1}]\cap R$ is a height one $\HH$-prime ideal of $R$ 
and by \cite[Proposition 2.9]{LLR} $R$ is an $\HH$-UFD. Therefore 
there exists a homogeneous prime element $v$ of $R$ such that 
\begin{equation}
\label{theta-equal}
\theta( u B) [x^{\pm 1}]\cap R = v R.
\end{equation}

It follows from \coref{nondeg0} that the degree of $v$ with respect to $x$ 
is at most 1. Thus exactly one of the following two cases holds: 
(1) $ v \in B$ or (2) the degree of $v$ with respect to $x$ equals $1$. 
First we prove several facts for the two cases simultaneously and 
then proceed with the rest separately. 

Denote the leading coefficient
of $v$ (as a left polynomial in $x$ with coefficients in $B$) by $u'$.
The first part of
\leref{new-prime-a} implies
\[
u B = \lc (\theta( u B) [x^{\pm 1}]\cap R) = \lc(v R) = u' B. 
\]
Therefore $u'$ is a homogeneous prime element of $B$ which is an 
associate of $u$ (in $B$).
Thus, after multiplying $v$ by a unit of $B$ we can assume that 
\begin{equation}
\label{uu'}
u'=u.
\end{equation}
Denote
\[
\theta(u) = u + c_{-1} x^{-1} + \cdots + c_{-m} x^{-m} 
\]
for some $m \in \Znn$, $c_{-1}, \ldots, c_{-m} \in B$, 
$c_{-m} \neq 0$. Let $i$ be the degree of $v$ as a polynomial 
in $x$ with coefficients in $B$ (i.e., $i=0$ in case (1) and $i=1$ 
in case (2)). It follows from \eqref{theta-equal} and the definition \eqref{C} of $\theta$
that 
\begin{equation}
\label{x1}
v = \sum_{l=n}^{i} \theta(u a_l) x^l 
\end{equation}
for some $n \in \Zset$ and $a_l \in B$ such that 
$n \leq i$, $a_i =1$, and $a_n \neq 0$.
Since 
$\theta(u) x^m \in \theta( u B) [x^{\pm 1}]\cap R$, equation \eqref{theta-equal}
and the fact that $R[x^{-1}] = B' [x^{\pm 1}; \sig]$ (where $B' = \theta(B)$) 
imply  
\begin{equation}
\label{x2}
\theta(u) x^m = v \Big( \sum_{l=m_1}^{m_2} \theta(b_l) x^l \Big) 
\end{equation}
for some $m_1 \le m_2$ in $\Zset$ and $b_l \in B$ such that 
$b_{m_1} \neq 0$ and $b_{m_2} \neq 0$. Let us 
substitute \eqref{x1} in \eqref{x2} and compare the 
coefficients of the powers of $x$, keeping in mind that 
$R[x^{-1}] = B' [x^{ \pm 1}; \sig]$. Since $B'$ is a 
domain, we obtain $n = i$, $m_1 = m_2$, $i + m_2 = m$, 
and 
\[
\theta(u) x^m = \theta(u) x^i \theta(b_{m_2}) x^{m_2}. 
\]
The $\HH$-equivariance of the map $\theta : B \to B'$ implies
$u = u \sig^i(b_{m_2})$. Thus, $b_{m_2}=1$ and
\begin{equation}
\label{vform} 
v = \theta(u) x^i.
\end{equation}

Case (1). Equations \eqref{uu'} and \eqref{vform} imply that in this case $v=u$. It follows 
from the definition \eqref{C} of the map $\theta$ and equation \eqref{vform} 
that $\de(u) = 0$. Hence,
\[
x u = \sig(u) x = \al u x.
\]

Case (2). By \eqref{uu'}, in this case $v = u x - c_\circ$ 
for some $c_\circ \in B$. Since $v$ is a normal element 
of $R = B[x;\sig, \de]$,
\[
v b = (\varphi_u \sig(b)) v, \; \; \forall b \in B.
\]
By a straightforward computation, this implies that $\de$ is given by 
\begin{equation}
\label{de9}
\de(b) = (u^{-1} c_\circ) b - \sig(b) (u^{-1} c_\circ), \; \; 
\forall b \in B.
\end{equation}
Invoking \eqref{vform} and using the definition \eqref{C} of the 
map $\theta$, we obtain
\[
u x - c_\circ = v = \theta(u) x = ux - (\la_\circ -1)^{-1} (\de \sig^{-1} (u)). 
\]
and $\de^2(u) =0$.
Hence,
\[
\de(u) = \al (\la_\circ -1) c_\circ \; \; 
\mbox{and} \; \; 
\de(c_\circ)=0.
\]
Again by a direct computation, one obtains from the 
first equation and the formula \eqref{de9} for 
$\de$ that 
\[
u c_\circ = (\al \la_\circ)^{-1} c_\circ u.
\]
Since $v$ is 
$\xh$-homogeneous, $\sig(u^{-1} c_\circ) = \la_\circ (u^{-1} c_\circ)$, and hence $\sig(c_\circ) = \al \la_\circ c_\circ$.   
It follows at once from \eqref{de9} that 
$\de(u^{-1} c_\circ) = - (\la_\circ-1) (u^{-1} c_\circ)^2$.
Using the formulas for $\de(u)$ and $\de(c_\circ)$, again 
by a direct computation one obtains
\[
x v = \al v x.
\]

The situation (ii) cannot occur for two different elements 
$c_\circ$ and $c'_\circ$, since in such a case 
$u x - c_\circ$ and $u x - c'_\circ$ would 
be two homogeneous prime elements of $R \setminus B$ which are 
not associates, thus contradicting \coref{nondeg0}. 

We have shown that case (1) implies situation (i) together with the stated additional conditions, and that case (2) implies situation (ii) together with its stated additional conditions. It remains to show that situations (i) and (ii) cannot occur simultaneously. This will follow from showing that $\de(u) = 0$ in situation (i) while $\de(u) \ne 0$ in situation (ii).

Assume situation (i), that is, $u$ is a prime element of $R$. The element $xu = \al ux + \de(u)$ must then lie in $uR$, whence $\de(u) \in uB$. This implies $\de(u) = 0$ by \leref{innerdelta-a}(a). 

Finally, assume situation (ii), that is, there is a homogeneous prime element $v = ux - c_\circ$ in $R$, for some $c_\circ \in B$. In particular, $c_\circ \ne 0$, and formula \eqref{pr-b} holds. Since $u^{-1}v = x - u^{-1}c_\circ$ is a homogeneous element of $R[E^{-1}]$, it must have the same $\xh$-degree $\la$ as $x$, and so $\sig(u^{-1}c_\circ) = \la_\circ u^{-1}c_\circ$. Hence, $\sig(c_\circ) = \al \la_\circ c_\circ$. Siince $vc_\circ = \al \la_\circ uc_\circ x - c_\circ^2$ lies in $Rv$, we must have $vc_\circ = c_\circ v$, whence $\al \la_\circ uc_\circ = c_\circ u$. Applying \eqref{pr-b} to $u$, we conclude that $\de(u) = \al (\la_\circ - 1) c_\circ \ne 0$.

This completes the proof of the theorem.
\end{proof}

\subsection{A classification of the homogeneous prime 
elements of $B[x;\sig,\de]$}
\label{2.5a}
By \cite[Proposition 2.9]{LLR} if the conditions 
{\rm(Cx1)--(Cx5)} are satisfied and $B$ is an $\HH$-UFD, 
then $R= B[x; \sig, \de]$ is an $\HH$-UFD. The following theorem 
classifies the homogeneous prime elements of $R$.

\bth{RtoBa}
Let $R=B[x; \sig, \de]$ be a skew polynomial algebra, equipped with a rational action of a $\KK$-torus $\HH$ by algebra automorphisms which leave $B$ invariant. Assume that $B$ is an $\HH$-UFD
and that $R$ satisfies {\rm(Cx1)--(Cx5)}.
Let $\{u_i \mid i \in I\}$ be a list of the homogeneous prime elements
of $B$ up to taking associates. Then there are the following three possibilities 
for a list of the homogeneous prime elements of $R$ up to taking associates:

{\rm (i)} $\{u_i \mid i \in I, i \neq i_0  \} \sqcup \{ v_{i_0}:=u_{i_0} x - c_\circ \}$
for some $i_0 \in I$ and $c_\circ \in B$ such that $u^{-1}_{i_0} c_\circ$ is a nonzero 
homogeneous element of $B[u^{-1}_{i_0}]$ with the same $\xh$-degree as $x$.

{\rm (ii)} $\{ u_i \mid i \in I \} \sqcup \{ x \}$.

{\rm (iii)} $\{ u_i \mid i \in I \}$.
\eth

We prove \thref{RtoBa} at the end of this subsection.
 
The next proposition describes the automorphisms corresponding to 
all homogeneous prime elements of $R$ in \thref{RtoBa} 
and the relationship of the skew derivation $\de$ to the 
homogeneous prime elements of $B$ and $R$. 
It is a direct consequence of \thref{R-B} and its proof is left to the reader.

\bpr{comm} Assume the setting of Theorem {\rm\ref{tRtoBa}}. For $i \in I$, let 
$\al_i \in \kx$ be such that $\sig(u_i) = h_\circ \cdot u_i = \al_i u_i$.
Then the following statements hold: 

{\rm (a)} In case {\rm (i)} of Theorem {\rm\ref{tRtoBa}}, we have 
\begin{gather*}
\de(b) = (u^{-1}_{i_0} c_\circ) b - \sig(b) (u^{-1}_{i_0} c_\circ), \; \forall b \in B,
\\
\de(u_i) = 0, \; \forall i \in I \setminus \{ i_0 \}, \; \; 
\de(u_{i_0}) = \al_{i_0}(\la_\circ-1) c_\circ \ne 0, \; \; 
\end{gather*}
and
\[
u_{i_0} c_\circ = (\al_{i_0} \la_\circ)^{-1} c_\circ u_{i_0}, \; \; 
\de(c_\circ) = 0, \; \; 
\de( u^{-1}_{i_0} c_\circ) = - (\la_\circ - 1) (u^{-1}_{i_0} c_\circ)^2.
\]
Furthermore,
\begin{align*}
&u_i x = \al_i^{-1} x u_i, \; \forall i \in I \setminus \{ i_0\}, \; \;   
v_{i_0} x = \al_{i_0}^{-1} x v_{i_0}, \; \; \mbox{and} \\
&v_{i_0} b =   (\varphi_{u_{i_0}} \sig(b)) v_{i_0} = \varphi_{u_{i_0}}( h_\circ \cdot b) v_{i_0}, \; 
\forall b \in B.
\end{align*}

{\rm (b)} In case {\rm (ii)} of Theorem {\rm\ref{tRtoBa}}, we have $\de =0$, 
\[
x b = \sig(b) x, \; \forall b \in B \; \; 
\mbox{and} \; \; u_i x = \al^{-1}_i x u_i, \; \forall i \in I.
\]

{\rm (c)} In case {\rm (iii)} of Theorem {\rm\ref{tRtoBa}}, we have $\de \neq 0$, 
\[
\de(u_i) = 0 \; \; 
\mbox{and} \; \; u_i x = \al^{-1}_i x u_i, \; \; 
\forall i \in I.
\]
\epr

\bre{nothird} In \thref{CGLind} we prove that the situation (iii) 
in \thref{RtoBa} can never occur when $B$ is a CGL 
extension. While we do not know an example when this 
situation can be realized, it appears to be difficult to rule it out 
in the generality of \thref{RtoBa}.
\ere

We will need the following lemma for the proof of \thref{RtoBa}.

\ble{deg0} If $u \in B$ is a prime element of $R = B[x; \sig, \de]$, 
then $u$ is a prime element of $B$.
\ele 

\begin{proof} The lemma follows from the fact that under the above assumptions, $u$ is normal in $B$ and 
$u B = (u R) \cap B$ is a completely prime ideal of $B$.
\end{proof}

\begin{proof}[Proof of Theorem {\rm\ref{tRtoBa}}]
By \thref{R-B}, for each $i \in I$, exactly one of the following 
two situations occurs:
\begin{enumerate}
\item[(*)] $u_i$ remains a prime element of $R$, or 
\item[(**)] there exists $c_\circ \in R$ such that $v_i:=u_i x - c_\circ$ is a 
homogeneous prime element of $R$. 
\end{enumerate}

\coref{nondeg0} implies that the situation (**) cannot occur for two different 
indices $i \in I$. Thus we have two cases: 
\begin{enumerate}
\item For all $i \in I$, (*) is 
satisfied.
\item For some $i= i_0 \in I$, (**) is satisfied, and 
for all $i \in I \setminus \{i_0\}$, (*) is satisfied.
\end{enumerate}
 Note that 
it is possible that $I = \varnothing$, in which case (1) holds.

Case (1). By \thref{R-B}, in this case $\de(u_i) = 0$ for all $i \in I$. 
It follows from \leref{deg0} that every homogeneous prime element
of $R$ of degree $0$ in $x$ is an associate of one of the elements 
$u_i$ for some $i \in I$. Again by \coref{nondeg0}, $R$ has 
no homogeneous prime elements of degree strictly greater than 1, and up to taking associates, 
$R$ has at most one homogeneous prime element of degree 1. Suppose that an element of the latter form exists, and denote it by $u x - c_\circ$ where $u$ is a homogeneous normal element of $B$. Set $d:= u^{-1}c_\circ$.
This is a homogeneous element of $B[E^{-1}]$ of degree $\la$.

If $u$ is a unit of $B$, then $x - d$ is a 
normal element of $R$. Hence, 
$\de$ is an inner $\sig$-der\-i\-va\-tion of $B$,
\begin{equation}
\label{some}
\de(b) = d  b - \sig(b) d, 
\; \; \forall b \in B.
\end{equation}
By \leref{innerdelta-a}(b), $d =0$. Furthermore, $\de =0$ 
and we are in the situation (ii).

If $u$ is not a unit of $B$, then after multiplying $u$ by a 
unit of $B$ we can assume that $u$ is equal to a product of
$u_i$'s (with at least one term). 
\thref{R-B} implies 
that in this case $\de(u_i) =0, \forall i \in I$ and 
thus $\de(u) = 0$. Moreover, $\de$ is still given by 
\eqref{some} and $\de(u) =0$ is equivalent to 
$u d = \al^{-1} d u$, 
where $\al \in \kx$ is such that $\sig(u) = \al u$.
It follows from \eqref{some} that $x-d$ is a normal 
element of $R[E^{-1}]$ and more precisely 
\[
(x-d) b = \sig(b) (x-d), \; \; \forall b \in B[E^{-1}].  
\]
Using this and the property $u d = \al^{-1} d u$, we obtain
\begin{align}
(u x - c_\circ) x &= u (x-d) \big[ (x-d) + d \big] = 
\big[ \al^{-1} (x-d) + \al^{-1} \la_\circ d  \big] u (x-d) 
\label{contr}
\\
&= \al^{-1} \big[ x + (\la_\circ -1) d \big] (u x - c_\circ).
\nn
\end{align}
This is a contradiction since $d \notin B$ and $\la_\circ \neq 1$.
Therefore in this subcase $B$ has no homogeneous prime elements of degree 
one in $x$ and we are in the situation (iii).

Case (2). \leref{deg0} implies that every homogeneous prime 
element of $R$ of degree $0$ is an associate of one of the 
elements $u_i$ for some $i \in I$. By \coref{nondeg0}, $R$ has no 
homogeneous prime elements of degree greater than 1 and 
all homogeneous prime elements of $R$ of degree 1 are 
associates of $u_{i_0} x - c_\circ$. It follows from \thref{R-B} 
that the element $u_{i_0} \in R$ is not prime since the two cases 
in that theorem are mutually exclusive.
Thus, in this case we 
are in the situation (i).
\end{proof}

\subsection{From $B$ to $B[x;\sig,\de]$ for strong $\HH$-UFDs}
\label{2.6a}
We complete this section with a treatment of the passage 
from $B$ to $R=B[x;\sig,\de]$ in the case when $B$ is a 
strong $\HH$-UFD. 

\bth{RtoBa-strong}
Let $R=B[x; \sig, \de]$ be a skew polynomial algebra, equipped with a rational 
action of a $\KK$-torus $\HH$ by algebra automorphisms which leave $B$ invariant. 
Assume that $B$ is a strong $\HH$-UFD and that $R$ satisfies {\rm(Cx1)--(Cx5)}.
Assume that if $\de =0$, then for every 
$\xi \in \kx$ there exists $t \in \HH$ such that $(t \cdot)|_B =\id_B$ and 
$t \cdot x = \xi x$. Then $R$ is a strong $\HH$-UFD.

Let $\{u_i \mid i \in I\}$ be a list of the $\HH$-prime elements
of $B$ up to taking associates, and let $t_i \in \HH$ be such that 
\[
u_i b = (t_i \cdot b) u_i, \; \; \forall b \in B. 
\]
Let $\al_i, \xi_i \in \kx$ be such that $\sig(u_i) = h_\circ \cdot u_i = \al_i u_i$
and $t_i \cdot x = \xi_i x$. Then there are three possibilities 
for a list of the $\HH$-prime elements of $R$ up to taking associates
exactly as in cases {\rm (i)--(iii)} of Theorem {\rm\ref{tRtoBa}}. 

In case {\rm (i)}, 
\begin{equation}
\label{strong1}
u_i c = (t_i \cdot c) u_i, \; \; 
\forall i \in I \setminus \{ i_0 \}, \; c \in R \quad 
\mbox{and} \quad v_{i_0} c = (t_{i_0} h_\circ \cdot c) v_{i_0}, 
\; \; \forall c \in R.
\end{equation}

In case {\rm (ii)}, for every $i \in I$ there exists $t'_i \in \HH$ such that 
$(t'_i \cdot)|_B =\id_B$ and $t'_i \cdot x = (\al_i \xi_i)^{-1} x$. There also
exists $t \in \HH$ such that $(t \cdot)|_B = \id_B$ and 
$t \cdot x = \la_\circ^{-1} x$. Furthermore,
\begin{equation}
\label{strong2}
u_i c = (t_i t'_i \cdot c) u_i, \; \; 
\forall i \in I, \; c \in R,
\end{equation}
and 
\begin{equation}
\label{strong3}
x c = (h_\circ t \cdot c) x,  \; \; \forall c \in R.
\end{equation}

In case {\rm (iii)}, 
\begin{equation}
\label{strong4}
u_i c = (t_i \cdot c) u_i, \; \; 
\forall i \in I, \; c \in R.
\end{equation}
\eth

\begin{proof} Since $B$ is an $\HH$-UFD, \thref{RtoBa} applies, and since $B$ is a strong $\HH$-UFD, we may assume that the list $\{ u_i \mid i\in I \}$ of $\HH$-prime elements of $B$ is also a complete list (up to taking associates) of the homogeneous prime elements of $B$. Once we show that each of the lists in the three cases of \thref{RtoBa} consists of $\HH$-normal elements of $R$, we will have established that $R$ is a strong $\HH$-UFD and will have verified the statement in the second paragraph of the theorem. Thus, all we need to prove is the validity 
of equations \eqref{strong1}--\eqref{strong4}. 

Equations \eqref{strong2} and \eqref{strong3}
are straightforward, in light of \prref{comm}(b), and are left to the reader. 

By \prref{comm}(a),
equation \eqref{strong1} is equivalent to showing
\[
\al_i \la(t_i) =1, \; \; \forall 
i \in I \setminus \{ 0 \} \quad
\mbox{and} \quad
\al_{i_0} \la_\circ \la(t_{i_0})=1
\]
in case (i). Using the fact that $u_{i_0}^{-1} c_\circ$ has $\xh$-degree
$\la$ (which follows from the homogeneity of $v_{i_0} = u_{i_0} x - c_\circ$) 
and \thref{R-B}, we obtain that for all $i \in I \setminus \{ i_0 \}$, 
\[
0 = \de(u_i) = (u_{i_0}^{-1} c_\circ) u_i - 
\al_i u_i (u_{i_0}^{-1} c_\circ) = 
(\la(t_i)^{-1} - \al_i) u_i u_{i_0}^{-1} c_\circ.
\]
Thus, $\al_i \la(t_i)=1$. Similarly, 
\[
\al_{i_0} (\la_\circ -1) c_\circ = 
\de(u_{i_0}) = 
(u_{i_0}^{-1} c_\circ) u_{i_0} - 
\al_{i_0} u_{i_0} (u_{i_0}^{-1} c_\circ) = 
(\la(t_{i_0})^{-1} - \al_{i_0}) c_\circ
\]
and $\al_{i_0} \la_\circ \la(t_{i_0}) =1$.

Finally, we consider case (iii). Applying \prref{comm}(c), 
we see that \eqref{strong4} is equivalent to
\begin{equation}
\label{aa}
\al_i \la(t_i) =1, \; \; \forall 
i \in I.
\end{equation}
Since $\de \neq 0$, there exists a homogeneous element $b$ of $B$ 
such that $\de(b) \neq 0$. Denote its degree by $\mu \in \xh$. Fix $i \in I$.
Acting by $\de$ on the equality $u_i b = \mu(t_i) b u_i$ and taking into 
account that $\de(u_i) = 0$ and $\sig(u_i) = \al_i u_i$ leads to 
\begin{equation}
\label{aa1}
\al_i u_i \de(b) = \mu(t_i) \de(b) u_i.
\end{equation}
On the other hand, $\de(b)$ is homogeneous of degree $\la \mu$ 
by \eqref{hdel}. Thus
\begin{equation}
\label{aa2}
u_i \de(b) = \la(t_i) \mu(t_i) \de(b) u_i.
\end{equation} 
Combining equations \eqref{aa1} and \eqref{aa2}, and using that 
$B$ is a domain and $\de(b) \neq 0$, leads to 
\[
\al_i^{-1} \mu(t_i) = \la(t_i) \mu(t_i).
\]
This implies \eqref{aa} and completes the proof of the theorem.
\end{proof} 

\sectionnew{Homogeneous prime elements and normal Gelfand--Tsetlin subalgebras of CGL extensions}
\label{CGL}
\subsection{Induction on homogeneous prime elements of CGL extensions}
\label{3.1}
Consider an iterated skew polynomial extension of length $N$,
\begin{equation} 
\label{itOre}
R := \KK[x_1][x_2; \sig_2, \delta_2] \cdots [x_N; \sig_N, \delta_N].
\end{equation}
For $k \in [0,N]$, denote the $k$-th algebra in the 
chain $R_k:=\KK[x_1][x_2; \sig_2, \de_2] 
\cdots [x_k; \sig_k, \de_k]$. Thus, $R_0 = \KK$ and 
$R_N = R$.

\bde{CGL} An iterated skew polynomial extension $R$ as in \eqref{itOre} 
is called a \emph{Cauchon--Goodearl--Letzter} (\emph{CGL}) \emph{extension} 
\cite[Definition 3.1]{LLR} if it is 
equipped with a rational action of a $\KK$-torus $\HH$ 
by $\KK$-algebra automorphisms satisfying the following conditions:
\begin{enumerate}
\item[(i)] For all $1 \leq j < k \leq N$, $\sig_k(x_j) = \la_{kj} x_j$
for some $\la_{kj} \in \kx$.
\item[(ii)] For every $k \in [2,N]$, $\de_k$ is a locally nilpotent 
$\sig_k$-derivation of $R_{k-1}$. 
\item[(iii)] The elements $x_1, \ldots, x_N$ are $\HH$-eigenvectors.
\item[(iv)] For every $k \in [1,N]$, there exists $h_k \in \HH$ such that 
$h_k \cdot x_k = \la_k x_k$ for some $\la_k \in \kx$, which is 
not a root of unity, and $h_k \cdot x_j = \la_{kj} x_j$, 
for all $j \in [1,k-1]$ (i.e., $\sig_k = (h_k \cdot)$ as elements 
of $\Aut(R_{k-1})$, for all $k \in [2,N]$). 
\end{enumerate}
\ede 

These conditions are chosen in such a way so the methods of \cite{GL} 
produce a finite stratification of $\Spec R$ by spectra 
of (commutative) Laurent polynomial rings and the 
deleting derivation method of \cite{Ca1} for studying $\Spec R$ 
is applicable.

For all CGL extensions $R$, the equality 
$\sig_k \delta_k =\la_k \delta_k \sig_k$ holds for
$k \in [2,N]$, just as in \eqref{sigdel}. Moreover,
$R_k = R_{k-1}[x_k; \sig_k, \de_k]$ is a 
Cauchon extension, for all $k \in [1,N]$: Conditions (Cx1)--(Cx4) are clear, and (Cx5) holds by \cite[Proposition 4.2]{GL}. Thus, all $\HH$-prime ideals of the algebras $R_1, \ldots, R_N$ are 
completely prime. 

The group of units of an iterated skew polynomial extension $R$
is reduced to scalars. Thus, two prime elements of $R$ 
are associates if and only if they are scalar multiples 
of each other.
By \cite[Proposition 3.2, Theorem 3.7]{LLR}, every CGL extension is an 
$\HH$-UFD, and every \emph{torsionfree} CGL extension is a UFD, where the latter property means that the subgroup of $\kx$ generated by $\{ \la_{kj} \mid 1\le j< k\le N \}$ is torsionfree. We address the strong $\HH$-UFD property in the following section.

The next theorem describes 
the homogeneous prime elements of $R=R_N$ iteratively 
from those of $R_{N-1}$. It proves that the situation (iii) 
from \thref{RtoBa} never arises in the framework of 
CGL extensions.

\bth{CGLind} Let $R$ be an arbitrary
CGL extension of length $N$ 
as in \eqref{itOre}. The following hold:

{\rm(a)} Let $\{u_i \mid i \in I\}$ be a list of the homogeneous prime elements of 
$R_{N-1}$ up to scalar multiples. 
There are two possibilities for a list of the homogeneous prime elements of 
$R$ up to scalar multiples:
\begin{enumerate}
\item[(i)] $\{u_i \mid i \in I \setminus \{i_0\} \} \sqcup \{ u_{i_0} x_N - c_\circ \}$
for some $i_0 \in I$ and $c_\circ \in R_{N-1}$ such that $u^{-1}_{i_0}c_\circ$ is a nonzero 
homogeneous element of $R_{N-1}[u^{-1}_{i_0}]$ 
with the same $X(\HH)$-degree as $x_N$.
\item[(ii)] $\{ u_i \mid i \in I \} \sqcup \{ x_N \}$. 
\end{enumerate}

{\rm(b)} Let $h \in \HH$, $\sig := (h \cdot) \in \Aut(R)$, and $\de$ be a locally nilpotent 
$\sig$-derivation of $R$ such that $\sig \de = q \de \sig$ 
for some $q \in \kx$ which is not a root of unity. If 
\begin{equation}
\label{de-assum}
\delta(u) = 0 \; \; 
\mbox{for all homogeneous prime elements $u$ of $R$},
\end{equation}
then $\de = 0$.
\eth

The proof of \thref{CGLind} will be given in subsection \ref{3.4}.

It follows from \thref{RtoBa} 
(or \thref{CGLind}) that a CGL extension $R$
has only a finite number of  pairwise nonproportional homogeneous 
prime elements. We will denote this number by $n$ and call it 
the {\em{rank}} of $R$. The rank of $R$ also equals the number of $\HH$-prime ideals of height $1$ in $R$.

For each $k \in [1,N]$, \thref{CGLind} in combination with \prref{comm} implies that $\rank R_k = \rank R_{k-1}$ if $\de_k \ne 0$, while $\rank R_k = \rank R_{k-1} + 1$ if $\de_k = 0$. Thus,
\begin{equation}
\label{rank}
\rank R_k = \big| \{ j \in [1,k] \mid \de_j = 0 \} \big|, \;\; \forall k \in [1,N].
\end{equation}

\subsection{Structure of CGL extensions}
\label{3.2}

Given a function $\mu : [1,N] \to [1,n]$, we 
define predecessor and successor functions
\[
p = p_\mu : [1,N] \to [1,N] \sqcup \{ - \infty \}, \qquad
s = s_\mu : [1,N] \to [1,N] \sqcup \{ + \infty \}
\]
for the level sets of $\mu$ by
\begin{equation}
\label{p}
p(k) = 
\begin{cases}
\max \{ j <k \mid \mu(j) = \mu(k) \}, 
&\mbox{if $\exists j < k$ such that $\mu(j) = \mu(k)$}, 
\\
p(k) = - \infty, \; & \mbox{otherwise} 
\end{cases}
\end{equation}
and
\begin{equation}
\label{s}
s(k) = 
\begin{cases}
\min \{ j > k \mid \mu(j) = \mu(k) \}, 
&\mbox{if $\exists j > k$ such that $\mu(j) = \mu(k)$}, 
\\
s(k) = + \infty, & \mbox{otherwise}. 
\end{cases}
\end{equation}

\bth{CGL} Let $R$ be a CGL extension of length $N$ and rank $n$ as 
in \eqref{itOre}. There exist a surjective function $\mu : [1,N] \to [1,n]$ 
and elements
\[
c_k \in R_{k-1} \; \; \mbox{for all} \; \; k \in [2,N] \; \; 
\mbox{with} \; \; p(k) \neq - \infty
\]
such that the elements $y_1, \ldots, y_N \in R$, recursively defined by 
\begin{equation}
\label{y}
y_k := 
\begin{cases}
y_{p(k)} x_k - c_k, &\mbox{if} \; \;  p(k) \neq - \infty \\
x_k, & \mbox{if} \; \; p(k) = - \infty,  
\end{cases}
\end{equation}
are homogeneous and have the property that for every $k \in [1,N]$,
\begin{equation}
\label{prime-elem}
\{y_j \mid j \in [1,k] , \; s(j) > k \}
\end{equation}
is a list of the homogeneous prime elements of $R_k$ up to scalar multiples.
\eth

\begin{proof} We define $\mu : [1,k] \twoheadrightarrow [1, \rank R_k]$ and elements $c_k \in R_{k-1}$ (when $p(k) \ne -\infty$) for $k= 1,\dots,N$. At each step, the new function $\mu$ will be an extension of the previous one, and so the corresponding new predecessor function $p$ will also be an extension of the previous one. However, the successor functions may change, so we will write $s_k$ for the successor function going with $\mu$ on $[1,k]$.

To start, set $\mu(1) := 1$. Note that $p(1) = -\infty$ and $s_1(1) = +\infty$. Moreover, $y_1 := x_1$ is the unique homogeneous prime element of $R_1$ up to scalar multiples.

Now let $1 < k \le N$, and assume that $\mu$ has been defined on $[1,k-1]$, together with elements $c_j \in R_{j-1}$ for $j\in [1,k-1]$ with $p(j) \ne -\infty$ and $y_j \in R_j$ for $j \in [1,k-1]$, such that the desired properties hold. In particular, $\{ y_j \mid j \in [1,k-1], \; s_{k-1}(j) \ge k \}$ is a list of the homogeneous prime elements of $R_{k-1}$ up to scalar multiples. There are two cases to consider, corresponding to cases (i), (ii) of \thref{CGLind}(a).

In the first case, there is some $j_0 \in [1,k-1]$ such that $s_{k-1}(j_0) \ge k$ and some $c_k \in R_{k-1}$ such that
\begin{equation}
\label{primeRk}
\{ y_j \mid j\in [1,k-1], \; j \ne j_0, \; s_{k-1}(j) \ge k \} \sqcup \{ y_{j_0} x_k - c_k \}
\end{equation}
is a list of the homogeneous prime elements of $R_k$ up to scalar multiples. In this case, $\rank R_k = \rank R_{k-1}$, and we extend $\mu$ to a function $[1,k] \twoheadrightarrow [1, \rank R_k]$ by setting $\mu(k) = j_0$. Since $s_{k-1}(j_0) \ge k$, we see that $p(k) = j_0$, and so $y_k :=  y_{j_0} x_k - c_k$. It is easily checked that the set \eqref{primeRk} equals \eqref{prime-elem}.

In the second case, $\{ y_j \mid j\in [1,k-1], \; s_{k-1}(j) \ge k \} \sqcup \{ x_k \}$ is a list of the homogeneous prime elements of $R_k$ up to scalar multiples. In this case, we set $\mu(k) = \rank R_k$ and readily check the desired properties.
\end{proof}

\thref{CGL} implies that the rank of each intermediate 
CGL extension $R_k$ is equal to 
\[
\big| \mu([1,k]) \big|.
\]   
The set \eqref{prime-elem} of pairwise nonproportional 
homogeneous prime elements of $R_k$ can be also written as
\[
\{ y_{ \max ( \mu^{-1}(i) \cap [1,k] ) } \mid 
i \in \mu([1,k]) \}.
\]

If $R$ is a strong $\HH$-UFD, then using \thref{RtoBa-strong} one easily describes 
the elements of $\HH$ which induce the automorphisms of each $R_k$ 
corresponding to the homogeneous prime elements \eqref{prime-elem} 
of $R_k$. We leave the details to the reader.

The following examples illustrate \thref{CGL}.

\bex{length3rank1}
Let $q \in \kx$ be a non-root of unity, and $r \in \Znn$. Let $R$ be the $\KK$-al\-ge\-bra given by generators $x_1$, $x_2$, $x_3$ and relations
\begin{align*}
x_2 x_1 &= q x_1 x_2 + 1,  &x_3 x_1 &= q^{-1} x_1 x_3 + x_2^r,  &x_3 x_2 &= q x_2 x_3.
\end{align*}
As noted in case 9 of \cite[Proposition 7.2.3.3]{Ric}, $R$ is an iterated skew polynomial ring
\[
R = \KK[x_1] [x_2; \sig_2,\de_2] [x_3; \sig_3,\de_3].
\]
There is a rational action of the rank one torus $\HH = \kx$ on $R$ by $\KK$-algebra automorphisms such that
\begin{align*}
\al.x_1 &= \al x_1,  &\al.x_2 &= \al^{-1} x_2,  &\al.x_3 &= \al^{-r-1} x_3
\end{align*}
for $\al \in \HH$. It is easily checked that $R$ is a  CGL extension. For later reference, we note that $R$ is not symmetric in the sense of \deref{symmetric}, since there is no element $h'_1 \in \HH$ such that $h'_1 \cdot x_2 = q^{-1} x_2$ and $h'_1 \cdot x_3 = qx_3$.

The elements $y_k$ of \thref{CGL} can be given as follows:
\begin{align*}
y_1 &= x_1,  &y_2 &= y_1 x_2 - (1-q)^{-1},  &y_3 &= y_2 x_3 - (q^{r+1} - 1)^{-1} x_2^{r+1}.
\end{align*}
For each $k=1,2,3$, the element $y_k$ is the unique homogeneous prime element of $R_k$ up to scalar multiples, and so $R_k$ has rank $1$. We leave the calculations to the reader.
\eex

\bex{quantSchubert} Let $\g$ be a finite dimensional simple Lie algebra 
and $\UU_q(\g)$ be the corresponding quantized universal enveloping algebra 
over a base field $\KK$ of arbitrary characteristic for a deformation 
parameter $q \in \KK^*$ which is not a root of unity. Denote the 
rank of $\g$ by $r$. Let $E_i$, $F_i$ and $K^{\pm1}_i$, $i \in [1,r]$ 
be the standard generators of $\UU_q(\g)$ as in \cite{BG}. Denote by 
$\UU^+$ and $\UU^-$ the 
subalgebras of $\UU_q(\g)$ generated by $E_1, \ldots, E_r$ and
$F_1, \ldots, F_r$, respectively. 
De Concini--Kac--Procesi \cite{DKP} and Lusztig \cite{L}
defined a family of subalgebras $\UU^\pm[w]$ of $\UU^\pm$ indexed 
by the elements $w$ of the Weyl group of $\g$ 
which in the case $\KK = \Cset$ are deformations 
of the coordinate rings of the corresponding Schubert cells
equipped with the standard Poisson structure.  
There is a canonical action 
of a $\KK$-torus $\HH$ of rank $r$ on $\UU_q(\g)$ that preserves 
all of the subalgebras $\UU^\pm[w]$. For each reduced expression 
of $w$, there is an iterated Ore extension presentation of $\UU^\pm[w]$ 
and this is a CGL extension presentation \cite{MC}. 
For later reference, we note that the part of the proof of \cite[Lemma 6.6]{Y-sqg} 
between equations (6.7) and (6.8) implies that 
$\UU^\pm[w]$ satisfy the hypotheses of \thref{full}. Because $\UU^\pm[w]$ satisfy the 
Levendorskii--Soibelman straightening law, these algebras 
are symmetric CGL extensions in the terminology of Section \ref{relation}.

The $\xh$-homogenous prime elements of $\UU^\pm[w]$ were 
explicitly described in \cite[Theorem 6.2(i)]{Y-sqg}. This theorem in 
particular implies that the rank of $\UU^\pm[w]$ equals the 
cardinality of the support of the Weyl group element $w$
(the set of 
simple roots $\al_i$ such that the corresponding reflection 
$s_i$ appears in one and thus in any reduced expression of $w$).
The recursive nature of the $\xh$-homogenous prime elements of $\UU^\pm[w]$ 
was proved in \cite[Proposition 3.3]{GeY} which
together with \cite[Theorem 6.2(i)]{Y-sqg} established the validity of \thref{CGL}
for the quantum Schubert cell algebras $\UU^\pm[w]$. The proofs of both results 
in \cite{GeY,Y-sqg} used 
in an essential way the second realization of the algebras 
$\UU^\pm[w]$ in terms of quantum function algebras from \cite[Theorem 3.7]{Y-lms} 
and Drinfeld's $R$-matrix commutation relations, while our proof 
of \thref{CGL} directly relies on the iterated Ore extension structure
of the algebra in discussion.
\eex

\thref{CGL} has important consequences for the 
structure of all CGL extensions $R$. 
In particular, it can be used to describe explicitly 
the normal subalgebra $\NN(R)$ and the normal Gelfand--Tsetlin
subalgebra $\NGZ(R)$ of $R$, cf. the introduction for definitions.
Recall that 
a matrix $\qb:=(q_{jk}) \in M_N(\kx)$ is called 
\emph{multiplicatively skewsymmetric} if 
\[
q_{jk} q_{kj} = q_{kk} = 1, \; \forall j,k \in [1,N].
\]
Such a matrix gives rise to the quantum affine space 
algebra 
\begin{equation}
\label{QAff}
\calA_{\qb} := \OO_{\qb}(\KK^N) = \frac{\KK \lcor Y_1, \ldots, Y_N \rcor}
{\lcor Y_k Y_j - q_{kj} Y_j Y_k \mid 1 \leq j < k \leq N \rcor}
\end{equation}
and the quantum torus 
\begin{equation}
\label{q-torus}
\calT_\qb := \OO_{\qb}((\kx)^N) = \calA_\qb[Y_1^{-1}, \ldots, Y_N^{-1}].
\end{equation}
Both $\calA_\qb$ and $\calT_\qb$ have Gelfand--Kirillov dimension 
equal to $N$. 

The algebra $R$ given by \eqref{itOre} has the $\KK$-basis 
\begin{equation}
\label{PBW}
\{ x^f := x_1^{m_1} \cdots x_N^{m_N} \mid 
f:= (m_1, \ldots, m_N) \in \Znn^N \}.
\end{equation}
Consider the reverse lexicographic order $\prec$ on $\Znn^N$: 
$(m'_1, \ldots, m'_N ) \prec (m_1, \ldots, m_N)$
if there exists $j \in [1,N]$ such that $m'_j < m_j$  
and $m'_k = m_k, \; \forall k \in [j+1, N]$. We will say that 
$b \in R \setminus \{ 0 \}$ 
has \emph{leading term} $\xi x^f$ where $\xi \in \kx$ and $f \in \Znn^N$ 
if 
\[
b = \xi x^f + \sum_{g \in \Znn^N,\; g \prec f} \xi_g x^g
\]
for some $\xi_g \in \KK$.
Denote $\lt(b) := \xi x^f$. It follows from condition (i) in \deref{CGL} 
that
\begin{equation}
\label{prod-lt}
\lt( x^f x^{f'} ) = 
\biggl( \prod_{k>j} \la_{kj}^{m_k m'_j} \biggr) x^{f + f'}, \; \; 
\forall \; f =(m_1, \ldots, m_N), 
f'=(m'_1, \ldots, m'_N) \in \Znn^N. 
\end{equation}
For $k \in [1,N]$, let $\ol{k} := (m_1, \ldots, m_N) \in \Znn^N$
be such that $m_j = 1$, if $j \leq k$ and 
$\mu(j) = \mu(k)$, and $m_j= 0$ otherwise.
Equation \eqref{y} implies
\begin{equation}
\label{tty}
\lt(y_k) = x^{\ol{k}}, \; \; \forall k \in [1,N].
\end{equation}

The elements $\la_{kj}$, $1 \leq j < k \leq N$ give rise to 
a unique multiplicatively skewsymmetric matrix 
\begin{equation}
\label{La-matrix}
\La = (\la_{jk}) \in M_N(\kx).
\end{equation}
Define the order function $O_- : [1,N] \to \Znn$ by
\[
O_-(k) := \max \{ m \in \Znn \mid p^m(k) \neq - \infty \},
\]
where as usual $p^0:= \id$. For $j,k \in [1,N]$ 
set
\begin{equation}
\label{al}
\al_{kj} := \prod_{m=0}^{O_-(j)} \la_{k, p^m(j)} \in \kx.
\end{equation}
It follows from \eqref{tty} and the homogeneity of $y_j$ (recall \deref{CGL}(i)) that 
\[
\sig_k \big( x^{\ol{j}} \big)= \al_{kj} x^{\ol{j}} \; \; 
\mbox{and} \; \; 
\sig_k(y_j) = \al_{kj} y_j, 
\; \; \forall \, 1 \leq j < k \leq N.
\]

Consider the multiplicatively skewsymmetric matrix $\qb = \qb(R) := (q_{jk}) \in M_N(\kx)$ 
such that
\begin{equation}
\label{q}
q_{jk}:= \prod_{i=0}^{O_-(j)} \prod_{m=0}^{O_-(k)} \la_{p^i(j), p^m(k)}
\end{equation}
for $j,k \in [1,N]$. It follows from \thref{CGL} and 
equations \eqref{prod-lt}, \eqref{tty} that 
\[
\lt(y_k y_j) = q_{kj} \lt(y_j y_k), \; \; \forall j,k \in [1,N].
\]
On the other hand, using \thref{CGL}, \prref{comm}, and the homogeneity of $y_j$, we see by induction on $k$ that $y_k$ and $y_j$ quasi-commute for $1\le j< k\le N$. Therefore,
\begin{equation}
\label{y_comrel}
y_k y_j = q_{kj} y_j y_k, \; \; \forall j,k \in [1,N]. 
\end{equation}
For $\qb = \qb(R)$, we thus obtain a $\KK$-algebra homomorphism
\[
\iota : \calA_\qb \to R \; \; 
\mbox{such that} \; \; \iota(Y_k) :=  y_k, \; \forall k \in [1,N],
\]
and in view of \eqref{tty} we see that $\iota$ is injective. It extends to an injective 
$\KK$-algebra homomorphism $\iota : \calT_\qb \to \Fract(R)$. 
We will identify 
$\iota(\calA_\qb)$ and $\iota(\calT_\qb)$ 
with $\calA_\qb$ and $\calT_\qb$. Denote by $\NN_\qb$ 
the subalgebra of $\calA_\qb$ generated by those $y_j$,
$j \in [1,N]$ such that $s(j) = + \infty$. 
It is obviously isomorphic to a quantum
affine space algebra of dimension $n$.

Recalling \prref{norm-UFD} and noting that the units of $R$ are scalars, we obtain the following result.

\bth{embed} Let $R$ be a CGL extension of length 
$N$ and rank $n$ as in Definition {\rm\ref{dCGL}}, and define $\qb = \qb(R)$ via {\rm\eqref{q}}.
The normal subalgebra $\NN(R)$ of
$R$ equals the quantum affine space algebra $\NN_\qb$ of dimension $n$.
The normal Gelfand--Tsetlin subalgebra $\NGZ(R)$ of $R$ 
corresponding to the canonical chain of subalgebras 
\[
R_1 \subset R_2 \subset \cdots \subset R_N
\]
equals the quantum affine space algebra $\calA_\qb$ of dimension $N$.
Furthermore, we have $\KK$-al\-ge\-bra embeddings
\[
\calA_\qb \subseteq R \subset \calT_\qb \subset \Fract(R),
\]
where $\calT_\qb$ is the quantum torus corresponding to $\calA_\qb$.
\eth

Applying \prref{comm} to the situation of \thref{CGL} leads to the 
following facts for the skew derivations $\de_k$ and commutation 
relations for the homogeneous prime elements of $R_k$ and 
the two terms of $y_k$.

\bpr{comm2} Keep the notation from Theorem {\rm\ref{tCGL}}. Let $k \in [2,N]$.

{\rm (a)} If $p(k) = - \infty$, then $\de_k =0$. 

The automorphisms of $R_k$ corresponding to its homogeneous prime elements are given by
\begin{align*}
y_j x_k &= \al_{kj}^{-1} x_k y_j, \; \forall j \in [1, k-1],  &y_k a &= \sig_k(a) y_k, \; \forall a \in R_{k-1},
\end{align*}
together with the actions of $\varphi_{y_j}$ on $R_{k-1}$ for $j\in [1,k-1]$ such that $s(j) \ge k$ {\rm(}obtainable by recursion{\rm)}.

{\rm (b)} If $p(k) \neq - \infty$, then the skew derivation 
$\de_k$ is nonzero and is given by
\[
\de_k(a) = y^{-1}_{p(k)} c_k a - \sig_k(a) y^{-1}_{p(k)} c_k, 
\; \forall a \in R_{k-1}
\]
and satisfies
\begin{align*}
\de_k(y_j) &= 0, \; \forall j \in [1,k-1] \; \; 
\mbox{such that} \; \; s(j) > k,   
&\de_k(y_{p(k)}) &= \al_{k p(k)}(\la_k-1) c_k \ne 0,  \\
\de_k(c_k) &= 0,   &\de_k( y_{p(k)}^{-1} c_k) &= - (\la_k -1) (y_{p(k)}^{-1} c_k)^2.
\end{align*}

The automorphisms of $R_k$ corresponding to its homogeneous 
prime elements are given by 
\begin{align*}
y_j x_k &= \al_{kj}^{-1} x_k y_j, \; \forall j \in [1, k-1] \; \; 
\mbox{such that} \; \; s(j) > k,  &y_k x_k &= \al_{k p(k)}^{-1} x_k y_k,  \\
y_k a &=   (\varphi_{y_{p(k)}} \sig_k(a)) y_k = \varphi_{y_{p(k)}}( h_k \cdot a) y_k, \; 
\forall a \in R_{k-1},
\end{align*}
together with the actions of $\varphi_{y_j}$ on $R_{k-1}$ for $j\in [1,k-1]$ such that $s(j) \ge k$ {\rm(}obtainable by recursion{\rm)}.
Furthermore, the components $y_{p(k)}$ and $c_k$ of $y_k$ satisfy
\[
y_{p(k)} c_k = (\al_{kp(k)} \la_k)^{-1} c_k y_{p(k)}.
\]
\epr

\bco{primeqcomm}
Every homogeneous prime element of $R$ quasi-commutes with $x_1,\dots, x_N$. More precisely,
\begin{equation}
\label{ycommx}
y_j x_k = \al_{kj}^{-1} x_k y_j \;\; \forall j,k \in [1,N] \text{\;with\;} s(j) > k.
\end{equation}
Consequently, $x_k \NN(R) = \NN(R) x_k$ for all $k \in [1,N]$.
\eco

\begin{proof} We just need to establish \eqref{ycommx}, since that implies the first statement, and the last statement follows from the first because $\NN(R)$ is generated by the homogeneous prime elements of $R$.

We proceed by induction on $l \in [1,N]$, to prove that \eqref{ycommx} holds for $j,k \in [1,l]$. The case $l=1$ is clear, since $y_1 = x_1$ and $\al_{11} = \la_{11} = 1$.

Now let $l>1$, and assume \eqref{ycommx} holds for $j,k \in [1,l-1]$. If $j \in [1,l-1]$ and $s(j) > l$, then both cases of \prref{comm2} yield $y_j x_l = \al_{lj}^{-1} x_l y_j$. Hence, it just remains to consider $y_l$.

If $p(l) = -\infty$, then $y_l = x_l$ and $\de_l = 0$. In this case, 
\[
y_l x_k = \sig_l(x_k) y_l = \la_{lk} x_k y_l = \al_{kl}^{-1} x_k y_l
\]
for $k \in [1,l-1]$, while $y_l x_l = x_l y_l = \al_{ll} x_l y_l$ because $\al_{ll} = \la_{ll} = 1$. Finally, suppose that $p(l) \ne -\infty$, and note that $\al_{ll} = \la_{ll}  \al_{l,p(l)} =  \al_{l,p(l)}$. Hence, it follows from \prref{comm2}(b) that $y_l x_l = \al_{ll}^{-1} x_l y_l$. Since $s(p(l)) = l$, our induction hypothesis implies that $y_{p(l)} x_k = \al_{k,p(l)}^{-1} x_k y_{p(l)}$ for $k \in [1,l-1]$. Appealing again to \prref{comm2}(b), we conclude that
\[
y_l x_k = (\varphi_{y_{p(l)}} \sig_l (x_k)) y_l = \al_{k,p(l)}^{-1} \la_{lk} x_k y_l = \al_{kl}^{-1} x_k y_l
\]
for $k \in [1,l-1]$.

This completes the induction.
\end{proof}

\subsection{Proof of \thref{CGLind}}
\label{3.4}
Given a positive integer $L$, we denote by \thref{CGLind}(a) $(L)$ the validity 
of the statement of \thref{CGLind} for all CGL extensions 
$R$ of length $N \leq L$. Similarly, we define 
the statements \thref{CGLind}(b) $(L)$, \thref{CGL} $(L)$, and \thref{embed} $(L)$.
The discussion above shows that
\begin{equation}
\label{implication0}
\mbox{ \thref{CGLind}(a) $(N)$ $\Longrightarrow$ Theorems \ref{tCGL} $(N)$ and \ref{tembed} $(N)$}.
\end{equation}
The proof of \thref{CGLind} will be completed once the following implications are established:
\begin{equation}
\label{implication1}
\mbox{ \thref{CGLind}(b) $(N-1)$ $\Longrightarrow $ \thref{CGLind}(a) $(N)$}.
\end{equation}
\begin{equation}
\label{last}
\mbox{\thref{CGLind}(a) $(N)$ $\Longrightarrow $ \thref{CGLind}(b) $(N)$}.
\end{equation} 

\begin{proof}[Proof of the implication \eqref{implication1}]
We apply \thref{RtoBa} to the skew polynomial algebra 
$R=R_{N-1}[x_N; \sig_N, \de_N]$. 
All we need to show is that in this setting, 
the situation (iii) in \thref{RtoBa} can never 
occur. Suppose that situation (iii) does occur. Then by \prref{comm}, $\de_N(u)= 0$ 
for all homogeneous prime elements $u$ of $R_{N-1}$ but 
$\de_N \neq 0$. However, this contradicts \thref{CGLind}(b) $(N-1)$.
\end{proof} 

Our proof of \eqref{last} involves some analysis of skew derivations on the quantum torus $\calT_\qb$. Given $f = (m_1,\dots,m_N) \in \Zset^N$, define the Laurent monomial 
\begin{equation}
\label{monom}
y^f := y_1^{m_1} \cdots y_N^{m_N} \in \calT_\qb.
\end{equation}
The algebra $\calT_\qb$ is $\Zset^N$-graded by 
\[
\deg y^f := f, \; \; \forall f \in \Zset^N. 
\]
We will say that a $\KK$-linear map $\eta : \calT_\qb \to \calT_\qb $ 
is \emph{$\Zset^N$-homogeneous of degree $g \in \Zset^N$} 
if $\eta(y^f) \in \KK y^{f+g}$ for all $f \in \Zset^N$.  
(The term homogeneous is already used in the context
of the $X(\HH)$-grading of $R$; we use the term 
$\Zset^N$-homogeneous to distinguish the two gradings.)
Given a general $\KK$-linear map $\eta : \calT_\qb \to \calT_\qb$, 
for $g \in \Zset^N$ there are uniquely defined $\Zset^N$-homogeneous 
$\KK$-linear maps $\eta^g$ of degree $g$ 
such that
\[
\eta = \sum_{g \in \Zset^N} \eta^g.
\] 
If $\sig$ is an automorphism of $\calT_\qb$ which preserves the $\Zset^N$-grading
and $\de$ is a $\sig$-derivation, then the component $\de^g$ is a $\sig$-derivation 
for each $g \in \Zset^N$. Since $\calT_\qb$ is finitely generated, $\de^g \neq 0$ 
for at most finitely many $g \in \Zset^N$.

Let $\prec$ be the reverse lexicographic order on $\Zset^N$ (defined as it was above on $\Znn^N$). Any nonzero element $u\in \calT_\qb$ can be uniquely written in the form
\begin{equation}  
\label{uyf}
u = \zeta_1 y^{f_1} +\cdots+ \zeta_r y^{f_r} \;\; \text{where\;} f_1 \prec \cdots \prec f_r \text{\;in\;}
 \Zset^N \text{\;and\;} \zeta_1, \dots, \zeta_r \in \kx.
 \end{equation}
We will say that $\zeta_r y^{f_r}$ is the \emph{leading term} of $u$ and denote it $\ltq(u)$, to distinguish it from our previous usage of leading terms. For future reference, observe that
\begin{equation}
\label{ltayb}
\ltq(a y_k + b) = \ltq(a) y_k \;\; \forall k \in [2,N], \; a,b \in \KK\langle y_1^{\pm1}, \dots, y_{k-1}^{\pm1} \rangle, \; a\ne 0.
\end{equation}

If $\sig$ and $\de$ are as above, with $\de \ne 0$, we have
\begin{equation}
\label{delcomp}
\de = \de^{g_1} +\cdots+ \de^{g_t} \;\; \text{with\;} g_1 \prec \cdots \prec g_t \text{\;in\;} \Zset^N, \text{\;all\;} \de^{g_i} \ne 0.
\end{equation}
Moreover, if $\sig \de = q \de \sig$ for some $q \in \kx$, then $\sig \de^{g_i} = q \de^{g_i} \sig$, for all $i\in [1,t]$.
For $m>0$ and $u$ as in \eqref{uyf}, the component of $\de^m(u)$ in degree $f_r+mg_t$ is $\zeta_r (\de^{g_t})^m (y^{f_r})$. Hence,
\begin{equation}
\label{del^m}
\de^m(u) = 0 \; \Longrightarrow \; (\de^{g_t})^m (\ltq(u)) = 0
\end{equation}
for all $m>0$ and nonzero $u\in \calT_\qb$.

\ble{delyyinv}
Let $\calT$ be a $\KK$-algebra domain, $\sig$ an automorphism of $\calT$, and $\de$ a $\sig$-de\-ri\-va\-tion on $\calT$ such that $\sig \de = q \de \sig$ for some $q \in \kx$ which is not a root of unity. Suppose $y \in \calT$ is a unit such that $\de^m(y) = \de^m(y^{-1}) = 0$ for some $m > 0$. Then 
$\de(y) = \de(y^{-1}) = 0$.
\ele

\begin{proof} We apply the idea of \cite[Lemme 7.2.3.2]{Ric}. Let $r,s \ge 0$ be maximal such that $\de^r(y) \ne 0$ and $\de^s(y^{-1}) \ne 0$. Then, by the $q$-Leibniz rule for $\de$ (e.g., \cite[Lemme 2.2]{Ca1}),
\begin{align*}
\de^{r+s}(yy^{-1}) &= \sum_{l=0}^{r+s} q^{l(l-r-s)} \tbinom{r+s}{l}_q \, \sig^{r+s-l} \de^l(y) \de^{r+s-l}(y^{-1})  \\
&= q^{-rs} \tbinom{r+s}{r}_q \, \sig^s \de^r(y) \de^s(y^{-1}) \ne 0.
\end{align*}
Since this is not possible when $r+s > 0$, we must have $r= s= 0$, and therefore $\de(y) = \de(y^{-1}) = 0$.
\end{proof}

\begin{proof}[Proof of the implication \eqref{last}]
Let $R$ be a CGL extension of length $N$ 
as in \deref{CGL}.
By \eqref{implication0}, the statements of Theorems \ref{tCGL}
and \ref{tembed} hold for $R$. 
Let $h \in \HH$ and $\sig = (h \cdot) \in \Aut(R)$.
Let $\de$ be a locally nilpotent $\sig$-derivation of $R$ satisfying
the conditions in \thref{CGLind}(b), and suppose that $\de \ne 0$.
The assumption \eqref{de-assum} implies
\begin{equation}
\label{delyk0}
\de(y_k) = 0, \; \; \forall k \in [1,N] \; \; 
\mbox{such that} \; \; s(k) = + \infty.
\end{equation} 

The automorphism $\sig = (h\cdot)$ of $R$ and the $\sig$-derivation $\de$ extend to an automorphism and $\sig$-derivation of $\Fract(R)$ satisfying $\sig \de = q \de \sig$.
 Obviously $\sig(\calT_\qb) = \calT_\qb$, and since $\de(y_k) \in R \subseteq \calT_\qb$ for all $k \in [1,N]$, we see that $\de(\calT_\qb) \subseteq \calT_\qb$. Now view $\sig$ and $\de$ as an automorphism and a $\sig$-derivation of $\calT_\qb$, note that $\sig$ preserves the $\Zset^N$-grading, and decompose $\de$ as in \eqref{delcomp}. Due to the assumption that $\de \ne 0$, we must have $\de(y_j) \ne 0$ for some $j \in [1,N]$. Moreover, $s(j) \ne +\infty$ and $\de^{g_t}(y_j) \ne 0$.
We will prove the following fact:
\begin{enumerate}
\item[(*)] {\em{If $\de^{g_t}(y_k)=0$ for some $k \in [1,N]$ with $p(k) \neq - \infty$, then 
$\de^{g_t}(y_{p(k)})=0$.}}
\end{enumerate}
By \eqref{delyk0}, $\de^{g_t} (y_k)=0$, for all $k \in [1,N]$ 
such that $s(k) = + \infty$.
A downward recursive application of (*) leads to $\de^{g_t}(y_k) =0$, for all $k \in [1,N]$, contradicting what we found above. This contradiction proves 
implication \eqref{last}.

We are left with showing (*). Assume that 
$\de^{g_t}(y_k) =0$ for some $k\in [1,N]$ with $p(k) \ne -\infty$. There exists $m>0$ such that 
\[
0 = \de^m(x_k) = \de^m( y^{-1}_{p(k)} y_k + y^{-1}_{p(k)} c_k).
\]
Since $c_k \in R_{k-1}$, the elements $c_k$ and $y_{p(k)}^{-1} c_k$ 
belong to the 
subalgebra of $\calT_\qb$ generated by $y_1^{\pm1}, \ldots, y_{k-1}^{\pm1}$. Observations  \eqref{ltayb} and \eqref{del^m} then imply $(\de^{g_t})^m (y^{-1}_{p(k)} y_k)=0$. 
But $\de^{g_t}(y_k)=0$, so
\begin{equation}
\label{der-00}
0 = (\de^{g_t})^m (y^{-1}_{p(k)} y_k) = \big((\de^{g_t})^m(y^{-1}_{p(k)}) \big) y_k, \; \;
\mbox{i.e.,} \; \; (\de^{g_t})^m(y^{-1}_{p(k)}) =0.  
\end{equation}
On the other hand, the restriction of $\de$ to $R$ is locally 
nilpotent. Hence, $\de^{m'}(y_{p(k)})=0$ for some $m'>0$. 
It follows from \eqref{del^m} that $(\de^{g_t})^{m'}(y_{p(k)})=0$.
We combine this and the second equality in \eqref{der-00}, and apply \leref{delyyinv}
to obtain $\de^{g_t}(y_{p(k)})=0$.
\end{proof}

\subsection{Freeness of $R$ over $\NN(R)$ for all CGL extensions $R$}
\label{3.5} 
Next we combine \thref{embed} with the method of the proof of 
\cite[Theorems 5.1 and 5.4]{Y-sqg} to obtain that 
every CGL extension $R$ is a free left and right module over
its normal subalgebra $\NN(R)$, with $\NN(R)$ as a 
direct summand.

Fix a CGL extension $R$ of length $N$ and rank $n$, and recall the 
definition of the surjective function $\mu : [1,N] \to [1,n]$ from \thref{CGL}. 
By \thref{embed}, $\NN(R)$ equals the quantum affine space algebra 
$\NN_\qb$ which is the subalgebra of $\calA_\qb$ generated by 
\[
\{ y_j \mid j \in [1,N], \; s(j) = + \infty \} = 
\{ y_{\max \; \mu^{-1}(i)} \mid i \in [1,n] \}.
\]
The Gelfand--Kirillov dimension of $\NN_\qb$ 
equals $n$.
Define the subset
\begin{equation}
\label{DeltaR}
\Delta(R) := \{ 
(m_1, \ldots, m_N) \in \Znn^N \mid 
\forall i \in [1,n] \; \exists k \in \mu^{-1} (i) 
\; \; \mbox{such that} \;\; m_k = 0 \},
\end{equation}
cf. \cite[Eq. (5.10)]{Y-sqg}.
Recall the definition of the vector $\ol{k} \in \Znn^N$ for 
$k \in [1,N]$ from Subsection \ref{3.2}. The following lemma is analogous 
to \cite[Lemma 5.3]{Y-sqg}. Its  
proof is left to the reader.

\ble{directsum} Let $R$ be a CGL extension of length $N$ and rank $n$, 
and $\mu : [1,N] \to [1,n]$ be the corresponding surjective function from 
Theorem {\rm\ref{tCGL}}. For every $f \in \Znn^N$, there exist unique 
$g \in \Delta(R)$ and $c_1, \ldots, c_n \in \Znn$ such that 
\begin{equation}
\label{basis}
f = g + c_1\, \ol{\max \; \mu^{-1}(1)} + \cdots + c_n\, \ol{ \max \; \mu^{-1}(n)}. 
\end{equation}
\ele

\bth{free} For all CGL extensions $R$, we have 
\[
R = \bigoplus_{g \in \Delta(R)} \NN(R) x^g = 
\bigoplus_{g \in \Delta(R)} x^g \NN(R), 
\]
cf. \eqref{PBW} and \eqref{DeltaR}.
\eth

Since $0 \in \Delta(R)$, \thref{free} has the following direct corollary.
The second statement follows from \coref{primeqcomm}, which implies that $\NN(R)x^g = x^g \NN(R)$ for all $g \in \Znn^N$.

\bco{fr} Every CGL extension $R$ is a free left and right module 
over its normal subalgebra $\NN(R)$, and $\NN(R)$ is an $(\NN(R),\NN(R))$-bimodule direct summand of $R$.
\eco  

\begin{proof}[Proof of Theorem \ref{tfree}]
 Recall the definition of the leading term of 
an element of $R$ from Subsection \ref{3.2}. The set 
\[
\{ x^f \mid f \in \Znn^N \} 
\]
is a $\KK$-basis of $R$. By \leref{directsum}, for every $f \in \Znn^N$ 
there exist unique $g \in \Delta(R)$ and $c_1, \ldots, c_n \in \Znn$ 
such that \eqref{basis} is satisfied. It follows from equations \eqref{prod-lt} 
and \eqref{tty} that 
\begin{align*}
\lt \big( y_{ \max \; \mu^{-1}(1)}^{c_1} \cdots y_{ \max \; \mu^{-1}(n)}^{c_n} x^g 
\big)
&= \xi x^f \; \; \mbox{and} 
\\
\lt \big( x^g y_{ \max \; \mu^{-1}(1)}^{c_1} \cdots y_{ \max \; \mu^{-1}(n)}^{c_n} 
\big)
&= \xi' x^f 
\end{align*}
for some $\xi, \xi' \in \kx$. This implies the statement of the theorem 
for filtration reasons analogously to the proof of \cite[Theorem 5.4]{Y-sqg}.
\end{proof}

\sectionnew{The strong $\HH$-UFD property}
\label{fullness}

We prove in this section that for any CGL extension $R$, the torus $\HH$ can be chosen so that $R$ is a strong $\HH$-UFD. That this does not hold for arbitrary choices of $\HH$ can be seen in the standard generic quantized coordinate ring of $\KK^N$, that is, the $\KK$-algebra $R$ with generators $x_1, \dots, x_N$ and relations $x_j x_k = q x_k x_j$ for $1\le j < k\le N$, where $q\in \kx$ is a non-root of unity. The rank $1$ torus $\HH = \kx$ acts rationally on $R$ with $\al \cdot x_j = \al x_j$ for all $\al \in \HH$ and $j \in [1,N]$. With this action, and with the iterated skew polynomial presentation
\[
R = \KK[x_1] [x_2;\sig_2] \cdots [x_N;\sig_N],
\]
$R$ is a torsionfree CGL extension. Thus, $R$ is an $\HH$-UFD and a UFD. However, it is not a strong $\HH$-UFD (assuming $N \ge 2$), because the homogeneous prime elements $x_j$ are not $\HH$-normal. This failure is easily repaired, however -- if $\HH$ is replaced by $(\kx)^N$ with its usual action, then $R$ becomes a strong $\HH$-UFD.

\subsection{Iteration of \thref{RtoBa-strong}}
The key to establishing the strong $\HH$-UFD property in a CGL extension is to have the hypothesis of \thref{RtoBa-strong} available at each step of the iteration.

\bth{full} Let $R$ be a CGL extension of length $N$ and rank $n$ as in Definition {\rm \ref{dCGL}}. Assume that for each 
$j \in [1,N]$ with $\de_j=0$ and each $\xi \in \kx$, there exists 
$t \in \HH$ such that $t \cdot x_j = \xi x_j$ and $t \cdot x_k = x_k$ for all $k \in [1,j-1]$. Then $R$ is a strong $\HH$-UFD.

The image of $\HH$ in $\Aut(R)$ is a $\KK$-torus of rank exactly $n$.
\eth

\begin{proof} It is an immediate corollary of  \thref{RtoBa-strong} that $R$ is a strong $\HH$-UFD. 

The kernel of the action map $\theta : \HH \rightarrow \Aut(R)$ is just the intersection of the kernels of the characters $\rho_j$ for $j \in [1,N]$, where $\rho_j \in \xh$ is the $\xh$-degree of $x_j$. Hence, $\ker \theta$ is a closed subgroup of $\HH$, and so $\theta(\HH)$ is a $\KK$-torus. This torus acts rationally on $R$, and $R$ is a CGL extension with respect to the $\theta(\HH)$-action. Moreover, $\theta(\HH)$ satisfies the hypotheses of the theorem. Thus, we may replace $\HH$ by $\theta(\HH)$, i.e., we may assume that the action of $\HH$ on $R$ is faithful.

Let
\[
\pi : \HH \longrightarrow \prod_{j\in D} \kx, \;\; \pi(h) = \bigl( \rho_j(h) \bigr)_{j\in D}
\]
be the natural projection, where $D := \{ j\in [1,N] \mid \de_j = 0 \}$. It is clear from the hypothesis of the theorem that $\pi$ is surjective. Since $n = |D|$ by \eqref{rank}, we thus have $\rank(\HH) \ge n$, and we will have equality once we establish that $\pi$ is injective.

From \eqref{hdel}, we have $(h\cdot)|_{R_{j-1}} \circ \de_j = \rho_j(h) \de_j \circ (h\cdot)|_{R_{j-1}}$ for $h\in \HH$ and $j \in [2,N]$. It follows that for $1\le k< j\le N$, either $\de_j(x_k) = 0$ or $\de_j(x_k)$ is homogeneous with $\xh$-degree $\rho_j + \rho_k$. If $\pi$ is not injective, there is a non-identity element $h\in \ker \pi$. Since $\HH$ acts faithfully on $R$, there must be some $j\in [1,N]$ such that $\rho_j(h) \ne 1$, and we may assume that $j$ is minimal for this property. Moreover, $j \notin D$ because $h \in \ker \pi$, so $j > 1$ and there is some $k \in [1,j-1]$ such that $\de_j(x_k) \ne 0$. Since $\de_j(x_k)$ is in $R_{j-1}$, its $\xh$-degree must be of the form $\sum_{i=1}^{j-1} m_i \rho_i$ for some $m_i \in \Znn$. On the other hand, this degree is $\rho_j+ \rho_k$ as noted above, so $\rho_j = -\rho_k + \sum_{i=1}^{j-1} m_i \rho_i$.  But the minimality of $j$ implies that $\rho_i(h) = 1$ for all $i < j$, and hence we obtain $\rho_j(h) = 1$, contradicting our assumption. Therefore $\pi$ is indeed injective.
\end{proof}

\subsection{Maximal tori and the strong $\HH$-UFD property}
We now describe the appropriate maximal torus for the strong $\HH$-UFD result. Let $R$ be a CGL extension of length $N$ as in \deref{CGL}. Equip $R$ with the rational action of the torus $(\kx)^N$ by invertible linear transformations given by the rule
\begin{equation}
\label{linact}
(\al_1, \dots, \al_N) \cdot (x_1^{m_1} \cdots x_N^{m_N}) = \al_1^{m_1} \cdots \al_N^{m_N} x_1^{m_1} \cdots x_N^{m_N}.
\end{equation}
The given action of $\HH$ on $R$ factors through the above $(\kx)^N$-action via a morphism (of algebraic groups) $\HH \rightarrow (\kx)^N$. Since nothing is lost by reducing $\HH$ modulo the kernel of its action, we may assume the action of $\HH$ on $R$ is faithful and then identify $\HH$ with its image in $(\kx)^N$. Thus, there is no loss of generality in assuming that $\HH$ is a closed subgroup of $(\kx)^N$. (For closedness, see, e.g., \cite[Corollary 1.4]{Bor}.)

Next, set
\begin{align}
\label{defG}
G &:= \{  \psi \in (\kx)^N \mid (\psi\cdot), \text{\;acting as in \eqref{linact}, is an automorphism of\;} R \}  \\
 &\,\cong \{ \zeta \in \Aut(R) \mid x_1,\dots,x_N \text{\;are\;} \zeta\text{-eigenvectors} \},  \notag
\end{align}
and observe that $G$ is a closed subgroup of $(\kx)^N$. Since $G$ is diagonalizable, its connected component of the identity, $G^\circ$, is a torus (e.g., \cite[Corollary 8.5]{Bor}). This subgroup is the unique maximal torus of $G$, and so it contains $\HH$. Let us set
\begin{equation}
\label{defHmax}
\Hmax(R) := G^\circ.
\end{equation}
(The definition of this group, and its position within $(\kx)^N$, depend on the given CGL extension presentation of $R$. However, we do not indicate this dependence in the notation.)
Since $\Hmax(R)$ contains $\HH$, the algebra $R$ is also a CGL extension with respect to $\Hmax(R)$. We shall see later that, in fact, $\Hmax(R) = G$ (i.e., $G$ is connected).

\begin{remark}
\label{universalHmax}
{\rm The group $\Hmax(R)$ associated with a CGL extension $R$ has the following universal property, assuming that we fix the CGL extension presentation \eqref{itOre} for $R$. If $\HH_1$ is any $\KK$-torus acting rationally on $R$ such that $(R,\HH_1)$ is CGL for the presentation \eqref{itOre}, then the action of $\HH_1$ on $R$ factors uniquely through the action of $\Hmax(R)$, via an algebraic group morphism $\HH_1 \rightarrow \Hmax(R)$. Thus, if we identify $\Hmax(R)$ with its natural image in $\Aut(R)$, the image of the action map $\HH_1 \rightarrow \Aut(R)$ must be contained in $\Hmax(R)$.}
\end{remark}

\bth{CGLfull} Let $R$ be a CGL extension of length $N$ and rank $n$ as in Definition {\rm \ref{dCGL}}, and assume that $\HH= \Hmax(R)$. Then $R$ is a strong $\HH$-UFD, and $\rank \HH = n$.
\eth

\begin{proof} Both conclusions will follow from \thref{full} once we verify the hypothesis of that theorem.

Let $j \in [1,N]$ such that $\de_j = 0$. By \prref{comm2}, $p(j) = -\infty$. There is some $l \in [1,N]$ with $s(l) = +\infty$ and $p^{O_-(l)}(l) = j$. By \coref{primeqcomm}, $\varphi_{y_l}(x_k) = \al_{kl}^{-1} x_k$ for all $k \in [1,N]$.

Set $\theta := \varphi_{y_l}^{-1} \circ \left( \prod_{m=0}^{O_-(l)} h_{p^m(l)} \cdot \right)$, which is an automorphism of $R$ for which all the $x_k$ are eigenvectors. Thus, we may identify $\theta$ with an element of the group $G \subseteq (\kx)^N$ considered above. Since $h_{p^m(l)}(x_k) = \la_{p^m(l),k} x_k$ when $k < p^m(l)$, we find that $\theta(x_k) = x_k$ for $k < j$ and $\theta(x_j) = \la_j x_j$. Hence, $\theta = (1,\dots,1,\la_j,*,\dots,*)$ as an element of $(\kx)^N$.

Let $\pi : (\kx)^N \rightarrow (\kx)^j$ be the projection onto the direct product of the first $j$ components of $(\kx)^N$. Then $\pi(\HH)$ is a closed subgroup of $(\kx)^j$. Set 
\[
J := \{ a \in (\kx)^j \mid a_k = 1\; \forall k < j\},
\]
another closed subgroup of $(\kx)^j$. Now $\pi(\HH) \cap J$ is a closed subgroup of $J$. Since $\HH$ has finite index in $G$, there is some $r>0$ such that $\theta^r \in \HH$. Then $\pi(\theta^{ri}) = (1,\dots,1,\la_j^{ri})$ lies in $\pi(\HH) \cap J$ for all $i \in \Zset$, and consequently $\pi(\HH) \cap J$ is infinite, because $\la_j$ is not a root of unity. However, $J$ is an irreducible $1$-dimensional variety, so we must have $\pi(\HH) \cap J = J$.

Thus, for any $\xi \in \kx$, there is some $t \in \HH$ with $\pi(t) = (1,\dots,1,\xi)$. Consequently, $t \cdot x_j = \xi x_j$ and $t \cdot x_k = x_k$ for all $k \in [1,j-1]$. This verifies the hypothesis of \thref{full}.
\end{proof}

The second conclusion of \thref{CGLfull} verifies a conjecture of Launois and Lenagan [unpublished].

\bco{normalHmax} 
Let $R$ be a CGL extension of length $N$ and rank $n$ as in Definition {\rm \ref{dCGL}}. For any nonzero normal element $u \in R$, there exists $h_u \in \Hmax(R)$ such that $ua = (h_u\cdot a)u$ for all $a \in R$.
\eco

\begin{proof}
We may assume that $\HH = \Hmax(R)$. Then $R$ is a strong $\HH$-UFD by \thref{CGLfull}, and the result follows from \thref{HnormalsinHUFD}.
\end{proof}

Now that we have established that the maximal torus associated to a CGL extension $R$ has rank equal to the rank of $R$, we can pin down the group $\Hmax(R)$ tightly, as follows.

\bth{pinHmax}
Let $R$ be a CGL extension of length $N$ and rank $n$ as in Definition {\rm \ref{dCGL}}, and define the group $G \subseteq (\kx)^N$ as in {\rm\eqref{defG}}. Then
\[
\Hmax(R) = G.
\]

Now set $D' := \{ k \in [2,N] \mid \de_k \ne 0 \} = \{ k \in [2, N] \mid p(k) \neq - \infty \}$. 
For each $k \in D'$, choose $j_k \in [1,k-1]$ such that $\de_k(x_{j_k}) \ne 0$, and choose $m_k \in \Znn^{k-1}$ such that the monomial $x^{m_k} = x_1^{m_{k1}} \cdots x_{k-1}^{m_{k,k-1}}$ appears in $\de_k(x_{j_k})$ when $\de_k(x_{j_k})$ is expressed in the PBW basis {\rm\eqref{PBW}}. Then
\begin{equation}
\label{tightHmax}
\Hmax(R) = \{ \psi = (\psi_1,\dots,\psi_N) \in (\kx)^N \mid \psi_k = \psi_{j_k}^{-1} \prod_{i=1}^{k-1} \psi_i^{m_{ki}}, \; \forall k \in D' \}.
\end{equation}
\eth

\begin{proof} By \thref{CGLfull} and equation \eqref{rank}, the rank of $\Hmax(R)$ is $n = N - |D'|$. Let $G_2$ denote the closed subgroup of $(\kx)^N$ described on the right hand side of \eqref{tightHmax}, and note that $G_2$ is a $\KK$-torus of rank $n$. We shall prove that $\Hmax(R) \subseteq G \subseteq G_2$. Since $\Hmax(R)$ and $G_2$ are connected groups of the same dimension, it will then follow that $\Hmax(R) = G_2$, proving both parts of the theorem.

By construction, $\Hmax(R) \subseteq G$, so only the inclusion $G \subseteq G_2$ remains. Let $\psi \in G$, and let $k \in D'$. On applying the automorphism $(\psi\cdot)$ to the relation 
\[
x_k x_{j_k} = \la_{k,j_k} x_{j_k} x_k + \de_k(x_{j_k}),
\]
we see that $\psi\cdot \de_k(x_{j_k}) = \psi_k \psi_{j_k} \de_k(x_{j_k})$. Consequently, all the monomials appearing in the PBW basis expansion of $\de_k(x_{j_k})$ must have $\psi$-eigenvalue $\psi_k \psi_{j_k}$. One of these is $x^{m_k}$, whose $\psi$-eigenvalue also equals $\prod_{i=1}^{k-1} \psi_i^{m_{ki}}$. Hence, $\psi_k = \psi_{j_k}^{-1} \prod_{i=1}^{k-1} \psi_i^{m_{ki}}$. This proves that $G \subseteq G_2$, as required.
\end{proof}

\sectionnew{Iterated sets of prime elements and the Cauchon quantum tori}
\label{relation} 
\subsection{Cauchon's quantum tori}
\label{4.1a}
Given a CGL extension $R$, there are now two completely 
different ways to embed it into a quantum torus. 
The first one is obtained via recursive applications 
of the Cauchon deleting derivation procedure \cite{Ca1}. 
The second is the one from \thref{embed} 
obtained via iterated sequences of homogeneous 
prime elements. These quantum tori are distinct
subalgebras of $\Fract(R)$. We relate them
for an important class of CGL extensions that 
contain all iterated skew polynomial extensions 
arising from quantized universal enveloping algebras 
of Kac--Moody algebras. More precisely,
if $R$ is a symmetric 
CGL extension as defined in \deref{symmetric}, then one 
can present $R$ as an iterated skew polynomial extension in 
two different ways by adjoining the variables in the orders 
$x_N, \ldots, x_1$ and $x_1, \ldots, x_N$. In \thref{rel}, we prove that 
the quantum tori obtained by applying the 
former procedure to the first iterated 
skew polynomial extension and 
the latter to the second one are equal as subalgebras 
of $\Fract(R)$. We furthermore derive an 
explicit formula expressing the iterated 
sets of homogeneous prime elements in terms 
of the Cauchon variables.

We begin by recalling the key steps in Cauchon's procedure 
of deleting derivations \cite[Section 3]{Ca1}. Let $R$ be a CGL extension 
of length $N$ as in \deref{CGL}. By abuse of notation, we will denote by the same symbol the extension of each $\sig_k$ to the automorphism $(h_k \cdot)$ of $\Fract(R)$, and the corresponding automorphism of any $\sig_k$-stable subalgebra of $\Fract(R)$.
For $l = N+1, \ldots, 2$,
Cauchon recursively defined
$N$-tuples of nonzero elements
\[
(x^{(l)}_1, \ldots, x^{(l)}_N)
\; \; 
\mbox{and subalgebras} \; \;  
R^{(l)} := \KK \lcor x^{(l)}_1, \ldots, x^{(l)}_N \rcor
\]
of $\Fract(R)$. In the first step $(l=N+1)$, one sets
\[
(x^{(N+1)}_1, \ldots, x^{(N+1)}_N) := 
(x_1, \ldots, x_N) \; \; \mbox{and} \; \; 
R^{(N+1)}=R.
\]

For the recursive step, let $l \in [2,N]$ and assume that $(x_1^{(l+1)}, \dots, x_N^{(l+1)})$ and $R^{(l+1)}$ have been defined. Cauchon proved that $R^{(l+1)}$ is an iterated skew polynomial ring of the form
\begin{equation}
\label{Rl-isom}
\begin{aligned}
R^{(l+1)} &= \KK[x_1^{(l+1)}] \cdots [x_l^{(l+1)}; \sig_l^{(l+1)}, \de_l^{(l+1)}] [x_{l+1}^{(l+1)}; \tau_{l+1}^{(l+1)}] \cdots [x_N^{(l+1)}; \tau_N^{(l+1)}]  \\
 &\cong \KK[x_1] [x_2; \sig_2,\de_2] \cdots [x_l; \sig_l, \de_l] [x_{l+1}; \tau_{l+1}] \cdots [x_N; \tau_N],
 \end{aligned}
\end{equation}
 where $\tau_k^{(l+1)} (x_j^{(l+1)}) = \la_{kj} x_j^{(l+1)}$ for $k > l,j$ and the isomorphism sends $x_j^{(l+1)} \mapsto x_j$ for all $j \in [1,N]$ \cite[Th\'eor\`eme 3.2.1]{Ca1}. In fact, as is easily seen, each $\sig_k^{(l+1)}$ and each $\tau_k^{(l+1)}$ is just $\sig_k$ acting on the appropriate subalgebra of $\Fract(R)$. To simplify the notation, one writes $\sig_j$, $\de_j$, $\tau_j$ for $\sig_j^{(l+1)}$, $\de_j^{(l+1)}$, $\tau_j^{(l+1)}$.
 The next $N$-tuples are defined by iterating the 
Cauchon map \eqref{C} in the following way:
\begin{equation}
\label{new-x}
x^{(l)}_j := 
\begin{cases}
x^{(l+1)}_j, 
& \mbox{if} \; \; j \geq l 
\\
\sum_{m=0}^\infty \frac{(1- \la_l)^{-m}}{[m]_{\la_l}!} 
\Big[ \delta_l^m \sig^{-m}_l \left(x^{(l+1)}_j \right) \Big]
\left(x^{(l+1)}_l \right)^{-m}, 
& \mbox{if} \; \; j < l.
\end{cases}
\end{equation}
In all cases the above sums are finite due 
to the local nilpotence of $\de_l$ and the commutation relation
$\sig_l \de_l = \la_l  \de_l \sig_l$. Cauchon [ibid] proved that 
\[
S_l := \Big\{ \left(x^{(l+1)}_l \right)^m 
\; \Big{|} \; m \in \Znn \Big\}
\]
is an Ore subset of $R^{(l)}$ and $R^{(l+1)}$ for $l \in [2,N]$ 
and that one has the following equality 
of $\KK$-subalgebras of $\Fract(R)$:
\begin{equation}
\label{isom-C}
R^{(l)}[S_l^{-1}] = R^{(l+1)} [S_l^{-1}].
\end{equation}

Denote the final $N$-tuple of Cauchon elements
\[
(\xbar_1, \ldots, \xbar_N) := (x^{(2)}_1, \ldots, x^{(2)}_N).
\]
By \eqref{isom-C}, we have the inclusions
\[
R \subset \KK \lcor \xbar_1^{\pm 1}, \ldots, \xbar_N^{\pm 1} \rcor
\subset \Fract (R).  
\]
Recall the definition \eqref{La-matrix} of the multiplicatively 
skewsymmetric matrix $\La = (\la_{jk}) \in M_N(\kx)$. 
It follows from \eqref{Rl-isom} that the algebra
$\KK \lcor \xbar_1^{\pm 1}, \ldots, \xbar_N^{\pm 1} \rcor$ is 
isomorphic to the quantum torus $\calT_\La$ by 
\begin{equation}
\label{qtorLa}
\xbar_k \mapsto Y_k, \; \; k \in [1,N]
\end{equation}
in terms of the notation \eqref{q-torus}. We note that only
very rarely is the quantum affine space algebra 
$\KK \lcor \xbar_1, \ldots, \xbar_N \rcor$ 
a subalgebra of $R$.

Repeated applications of the $\HH$-equivariance of the Cauchon map \cite[Lemma 2.6]{LLR} give the following:
\begin{equation}
\label{xjbarH}
\xbar_j \text{\;is an\;} \HH\text{-eigenvector with the same\;} \HH\text{-eigenvalue as\;} x_j, \;\; \forall j \in [1,N].
\end{equation}
(We do not call $\xbar_j$ homogeneous here, because $\Fract(R)$ is not $\xh$-graded.)

For $j,k \in [1,N]$, let $R_{[j,k]}$ be the unital subalgebra of $R$ generated 
by those $x_i$ for which $j \leq i \leq k$. In other words
\[
R_{[j,k]}:= \KK \lcor x_j, \ldots, x_k \rcor \; \; 
\mbox{if} \; \; j\leq k \; \; \mbox{and} \; \; 
R_{[j,k]}:= \KK \; \; \mbox{otherwise}.  
\]

\ble{xjbarxj} Suppose that $\de_k(x_j) \in R_{[j+1,k-1]}$ for $1\le j <k \le N$. Then $\xbar_j - x_j \in \Fract R_{[j+1,N]}$ for all $j \in [1,N]$.
\ele

\begin{proof}
For $1\le j< l\le N$, it follows from \eqref{new-x} and our hypothesis that
\[
x_j^{(l)} - x_j^{(l+1)} \in \Fract \KK \langle x_{j+1}^{(l+1)}, \dots, x_l^{(l+1)} \rangle.
\]
A downward induction on $m$ then yields
\[
\Fract \KK \langle x_k^{(m)}, \dots, x_N^{(m)} \rangle \subseteq \Fract R_{[k,N]}
\]
for $k,m \in [1,N]$, whence $x_j^{(l)} - x_j^{(l+1)} \in \Fract R_{[j+1,N]}$ for $1\le j< l\le N$. Since $x_j^{(N+1)} = x_j$ and $x_j^{(j+1)} = \xbar_j$, the lemma follows.
\end{proof}

\subsection{Symmetric CGL extensions}
\label{4.2a}
For an iterated skew polynomial extension $R$ as in \eqref{itOre},
the notation above connects with the notation from Subsection \ref{3.1} as follows: $R_k = R_{[1,k]}$ for $k \in [1,N]$. 
Set
\[
R'_k := R_{[k, N]}.
\] 

\bde{symmetric} We call a CGL extension $R$ of length $N$ as in 
\deref{CGL} {\em{symmetric}} if the following two conditions hold:
\begin{enumerate}
\item[(i)] For all $1 \leq j < k \leq N$,
\[
\de_k(x_j) \in R_{[j+1, k-1]}.
\]
\item[(ii)] For all $j \in [1,N]$, there exists $h'_j \in \HH$ 
such that 
\[
h'_j \cdot x_k = \la_{kj}^{-1} x_k = \la_{jk} x_k, \; \forall 
k \in [j+1, N]
\]
and $h'_j \cdot x_j = \la'_j x_j$ for some $\la'_j \in \kx$ which is not a root of unity.
\end{enumerate}
\ede

Given a symmetric CGL extension $R$ as in \deref{symmetric}, both $\sig_k$ and $\de_k$ preserve $R_{[j,k-1]}$ for 
$j \in [1,k-1]$. Now let $j \in [1, N-1]$, fix a choice of $h'_j \in \HH$ as in  the definition, and denote 
\[
\sig'_j := (h'_j \cdot)|_{R'_{j+1}}.
\]
This is an automorphism of $R'_{j+1}$, uniquely determined by the conditions
\[
\sig'_j(x_k) = \la_{jk} x_k, \; \; \forall k \in [j+1,N],
\]
even though $h'_j$ may not be unique. The conditions (i) and (ii) imply that the inner 
$(h'_j \cdot)$-derivation on $R$ given by $a \mapsto x_j a - (h'_j \cdot a) x_j$ restricts to a $\sig'_j$-derivation $\de'_j$ of $R'_{j+1}$ such that 
\[
\de'_j(x_k): = x_j x_k - \la_{jk} x_k x_j = - \la_{jk} \de_k(x_j), \; \; 
\forall k \in [j+1, N].
\] 
The maps $\sig'_j$ and 
$\de'_j$ preserve $R_{[j+1,k]}$ for $k \in [j+1, N]$.

The following lemma and corollary are straightforward.

\ble{two-Ore} For all symmetric CGL extensions $R$ of length 
$N$ as in Definition {\rm\ref{dsymmetric}}, we have the skew polynomial extensions
\[
R_{[j,k]} = R_{[j,k-1]}[x_k; \sig_k, \de_k] \; \; 
\mbox{and} \; \; 
R_{[j,k]} = R_{[j+1,k]} [ x_j; \sig'_j, \de'_j]
\]
for all $1 \leq j < k \leq N$. Both derivations $\de_k$ and $\de'_j$ are locally nilpotent.
\ele

\bco{secondpres} Every symmetric CGL extension $R$ 
as above 
has the second CGL extension presentation 
\begin{equation}
\label{itOre2}
R = \KK [x_N] [x_{N-1}; \sig'_{N-1}, \de'_{N-1}] 
\cdots [x_1; \sig'_1, \de'_1]
\end{equation}
by keeping the same $\KK$-torus $\HH$ as in the original presentation 
and using the elements $h'_N, \ldots, h'_1 \in \HH$.
\eco 

\bre{moreOre} More generally, a symmetric CGL extension $R$ 
as above has a CGL extension presentation associated to every 
permutation $\tau$ of $\{1, \ldots, N \}$ 
such that 
\[
\tau(k) = \max \, \tau( [1,k-1]) +1 \; \;
\mbox{or} \; \; \tau(k) = \min \, \tau( [1,k-1]) - 1, 
\; \; \forall k \in [2,N].
\] 
This presentation is given by
\[
R = \KK [x_{\tau(1)}] [x_{\tau(2)}; \sig''_{\tau(2)}, \de''_{\tau(2)}] 
\cdots [x_{\tau(N)}; \sig''_{\tau(N)}, \de''_{\tau(N)}],
\]
where $\sig''_{\tau(k)} := \sig_{\tau(k)}$ and 
$\de''_{\tau(k)} := \de_{\tau(k)}$ if 
$\tau(k) = \max \, \tau( [1,k-1]) +1$, while 
$\sig''_{\tau(k)} := \sig'_{\tau(k)}$ and 
$\de''_{\tau(k)} := \de'_{\tau(k)}$ if 
$\tau(k) = \min \, \tau( [1,k-1]) -1$.
For this presentation, one keeps the original $\KK$-torus
$\HH$ and uses the elements $h''_{\tau(1)}, \ldots, h''_{\tau(N)} \in \HH$
given by $h''_{\tau(k)} = h_{\tau(k)}$ if $\tau(k) = \max \, \tau( [1,k-1]) +1$  
and $h''_{\tau(k)} = h'_{\tau(k)}$ if $\tau(k) = \min \, \tau( [1,k-1]) -1$.
The original CGL extension presentation of $R$ corresponds to $\tau=\id$, and 
the one in \coref{secondpres} to the permutation given by $\tau(k) = N+1 -k$, 
for all $k \in [1,N]$.
\ere

\subsection{The relationship theorem}
\label{4.3a}
Let $R$ be a CGL extension of length $N$ as in \deref{CGL}. 
Equation \eqref{y} defines the iterated set of elements $y_1, \ldots, y_N$ of $R$, which incorporate the homogeneous prime elements of the intermediate 
algebras $R_k$ by \eqref{prime-elem}. Recall from Subsection \ref{3.2} that the subalgebra 
$\KK\lcor y_1^{\pm 1}, \ldots, y_N^{\pm 1} \rcor$ of $\Fract(R)$
is a quantum torus denoted by $\calT_{\qb}$ where the multiplicatively 
skewsymmetric matrix $\qb = (q_{jk}) \in M_N(\kx)$ is given by \eqref{q}. 

The Cauchon procedure applied to the CGL extension presentation 
\eqref{itOre} of $R$ produces the elements $\xbar_1, \ldots, \xbar_N \in \Fract(R)$. 
The subalgebra $\KK \lcor \xbar^{\pm 1}_1, \ldots, \xbar_N^{\pm 1} \rcor $ 
of $\Fract (R)$ is isomorphic to the quantum torus $\calT_\La$ associated to 
the multiplicatively skewsymmetric matrix $\La = (\la_{jk}) \in M_N(\kx)$, see
\eqref{La-matrix}. The quantum tori 
$\calT_\qb= \KK\lcor y_1^{\pm 1}, \ldots, y_N^{\pm 1} \rcor$ 
and $\calT_\La = \KK \lcor \xbar^{\pm 1}_1, \ldots, \xbar_N^{\pm 1} \rcor $ 
are completely different subalgebras of $\Fract(R)$ except 
in some very special cases. For instance, in \exref{length3rank1} one calculates that
\begin{align*}
\xbar_1 &= x_1 + q (1-q^{r+1})^{-1} x_2^r x_3^{-1} + (q-1)^{-1},  &\xbar_2 &= x_2,  &\xbar_3 &= x_3.
\end{align*} 
Moreover, $y_3 = \xbar_1 \xbar_2 \xbar_3$ but $y_2 = \xbar_1 \xbar_2 - (1-q^{r+1})^{-1} \xbar_2^{r+1} \xbar_3^{-1}$, which is a non-unit in $\calT_\La$. Thus, in this example we have $\calA_\qb \subseteq \calT_\La$ but $\calT_\qb \nsubseteq \calT_\La$.

In order to relate the two embeddings of a CGL extension in
quantum tori, we need to use a reverse CGL extension 
presentation in one of the two cases. For this 
reason, from now on we assume that $R$ is a 
symmetric CGL extension. Then we have the 
second CGL extension presentation of $R$ from 
\eqref{itOre2}. We apply the Cauchon procedure 
to this CGL extension presentation, keeping the indices in descending order: $N,N-1,\dots, 1$. The analogs 
of the elements $\xbar_1, \ldots, \xbar_N$ 
for this case will be denoted by 
\[
\xbar'_N, \ldots, \xbar'_1 \in \Fract(R).
\]
The discussion in Subsection \ref{4.1a} shows that 
\[
\calT'_\La := \KK \lcor (\xbar'_N)^{\pm 1}, \ldots, (\xbar'_1)^{\pm 1} \rcor 
\subset \Fract(R) 
\]
is a quantum torus which is isomorphic to $\calT_\La$. 
In particular, 
\[
\xbar'_j \xbar'_k = \la_{jk} \xbar'_k \xbar'_j, 
\; \; \forall j, k \in [1,N].
\] 

Since the presentation \eqref{itOre2} satisfies $\de'_j(x_k) \in R_{[j+1,k-1]}$ for $1\le j< k\le N$, observation \eqref{xjbarH} and \leref{xjbarxj} yield the following for all $k \in [1,N]$:
\begin{equation}
\label{xkbar'}
\begin{aligned}
 &\xbar'_k \text{\;is an\;} \HH\text{-eigenvector with the same\;} \HH\text{-eigenvalue as\;} x_k.  \\
 &\xbar'_k - x_k \in \Fract R_{[1,k-1]}.
\end{aligned}
\end{equation}

The main result of this section relates the embeddings 
$R \subset \calT_\qb$ and $R \subset \calT'_\La$, and the 
generators $y_1, \ldots, y_N$ of $\calT_\qb$ and 
$\xbar'_1, \ldots, \xbar'_N$ of $\calT'_\La$.
From another perspective it gives explicit simple formulas for 
the homogeneous prime elements of the intermediate 
algebras $R_k$ in terms of the Cauchon elements for $R$.
The proof of the theorem relies on an extension 
of the argument of the proof of \cite[Theorem 3.1]{GeY} 
of Geiger and Yakimov.

\bth{rel} Let $R$ be a symmetric CGL extension of length $N$
as in Definition {\rm\ref{dsymmetric}}. Then 
\begin{equation}
\label{equal-tor}
\calT'_\La = \calT_\qb
\end{equation}
and for all $k \in [1,N]$,
\begin{equation}
\label{ytox}
y_k = \xbar'_{p^{O_-(k)}(k)} \cdots \xbar'_k.
\end{equation}
\eth

\begin{proof} Obviously \eqref{equal-tor} follows from \eqref{ytox}.
We prove \eqref{ytox} by induction on $k$. For $k=1$, we have 
$\xbar'_1 = x_1$ by the current Cauchon procedure, and 
$y_1 = x_1$ by \eqref{y}.

Assume \eqref{ytox} holds for all $k = 1, \ldots, l-1$
for some $l \in [2,N]$. Therefore 
\[
\xbar'_k = 
\begin{cases}
y_{p(k)}^{-1} y_k, & \mbox{if} \; \; p(k) \neq - \infty
\\
y_k, & \mbox{if} \; \; p(k) = - \infty
\end{cases}
\]
for all $k \in [1,l-1]$. Since $y_l y_k = q_{lk} y_k y_l$, it 
follows from \eqref{q} that 
\[
y_l \xbar'_k = \al_{kl}^{-1} \xbar'_k y_l, \; \; \forall k \in [1,l-1],
\]
in terms of the scalars $\al_{kl}$ defined 
in \eqref{al}. At the same time, working in the quantum torus 
$\calT'_\La$ we obtain
\[
\big( \xbar'_{p^{O_-(l)}(l)} \cdots \xbar'_l \big) \xbar'_k = 
\al_{kl}^{-1} \xbar'_k \big( \xbar'_{p^{O_-(l)}(l)} \cdots \xbar'_l \big)
\]
for the same values of $k$. Set 
\[
z_l := \big( \xbar'_{p^{O_-(l)}(l)} \cdots \xbar'_l \big)^{-1} y_l = (\xbar'_l)^{-1} y_{p(l)}^{-1} y_l \in \Fract(R),
\]
taking account of \eqref{ytox} for $k = p(l)$ and \eqref{xkbar'}.
We have
\begin{equation}
\label{cent1}
z_l \xbar'_k = \xbar'_k z_l, \; \; \forall k \in [1,l-1].
\end{equation}
Applying again \eqref{xkbar'} and \eqref{ytox} gives
\begin{equation}
\label{FracRl-1}
\Fract(R_{l-1}) = \Fract( \KK \langle y_1, \dots, y_{l-1} \rangle) = \Fract( \KK  \langle \xbar'_1, \dots, \xbar'_{l-1} \rangle).
\end{equation}
Consequently, \eqref{cent1} implies
\begin{equation}
\label{zlcomm}
z_l b = b z_l, \; \; \forall b \in \Fract(R_{l-1}).
\end{equation}

By \eqref{y} and \eqref{xkbar'},
\begin{equation}
\label{frac-eqs}
z_l \in \Fract(R_l).
\end{equation}
Equation \eqref{y} together with \eqref{xkbar'} 
implies 
\[
z_l = 1 +(\xbar'_l)^{-1} b_l 
\]
for some $b_l \in \Fract(R_{l-1})$.
Now it follows from \eqref{zlcomm} and the fact that $z_l$ commutes with 
itself that $z_l \xbar'_l = \xbar'_l z_l$, and hence
\[
z_l x_l = x_l z_l,
\]
due to \eqref{xkbar'}.
Applying again \eqref{zlcomm} and \eqref{frac-eqs} gives
$z_l \in \ZZ(\Fract(R_l))$. The homogeneity of the right hand side of 
\eqref{y} and the fact that $\xbar'_k $ has the same $\HH$-eigenvalue 
as $x_k$, for all $k \in [1,N]$, imply $z_l$ is fixed under the action of $\HH$, that is,
\[
z_l \in \ZZ(\Fract(R_l))^\HH.
\]
But $R_l$ is a CGL extension and so the strong $\HH$-rationality 
result \cite[Theorem II.6.4]{BG} implies that 
$\ZZ(\Fract(R_l))^\HH = \KK$. This forces $z_l = 1$, which establishes
the validity of \eqref{ytox} for $k=l$.
\end{proof}

\section*{Acknowledgement}

We thank the anonymous referees for their careful reading of the manuscript and helpful suggestions.


\end{document}